\documentclass[MSNbibl,number,citesort,dvips]{arxbj}
\usepackage{subenv}
\usepackage{mathbh,mathrsfs}

\aid{0}
\volume{18}
\issue{4}
\pubyear{2012}
\firstpage{1361}
\lastpage{1385}
\doi{10.3150/11-BEJ369}

\makeatletter

\define@key{bid}{pmid}{}

\newtheorem{theorem}{Theorem}
\newtheorem{result}{Result}
\newtheorem{lemma}{Lemma}
\newproclaim{definition}{Definition}
\newtheorem{assumption}{A\hspace*{-2pt}}
\newremark{remark}{Remark}
\newtheorem{property}{Property}
\newtheorem{corollary}{Corollary}
\def\VV#1{\operatorname{Var}}
\makeatother

\begin{document}
\begin{frontmatter}

\title{Uniform convergence of the empirical cumulative
distribution function under informative selection from a finite population}
\runtitle{Uniform c.d.f. convergence under informative selection}
%
%
%
\begin{aug}
\author[ad1]{\fnms{Daniel} \snm{Bonn\'ery}\thanksref{ad1,e1,u1}\ead[label=e1,mark]{daniel.bonnery@ensai.fr}\ead[label=u1,url,mark]{www.ensai.com/daniel-bonnery-rub,65.html}},
\author[ad2]{\fnms{F. Jay} \snm{Breidt}\thanksref{ad2}\ead[label=e2]{jbreidt@stat.colostate.edu}\ead[label=u2,url]{www.stat.colostate.edu/\textasciitilde jbreidt}}
\and
\author[ad1]{\fnms{Fran\c{c}ois} \snm{Coquet}\thanksref{ad1,e3,u3}\ead[label=e3,mark]{fcoquet@ensai.fr}\ead[label=u3,url,mark]{www.ensai.com/francois-coquet-rub,30.html}}

\runauthor{D. Bonn\'ery, F.J. Breidt and F. Coquet}

\address[ad1]{Ensai, Campus de Ker-Lann, Rue Blaise Pascal -- BP 37203,
35172 Bruz -- cedex,
France.\\
\printead{e1,e3},\\ \printead{u1}, \printead*{u3}}

\address[ad2]{Department of Statistics,
Colorado State University,
Fort Collins, CO 80523-1877, USA. \\
\printead{e2}, \printead{u2}}
\end{aug}

\received{\smonth{10} \syear{2010}}
\revised{\smonth{3} \syear{2011}}

%
\begin{abstract}
Consider informative selection of a sample from a finite population.
Responses are realized as independent and
identically distributed (i.i.d.) random variables with a probability
density function (p.d.f.) $f$, referred to
as the superpopulation model. The selection is informative in the sense
that the sample responses,
given that they were selected, are not i.i.d. $f$. In general, the
informative selection mechanism may induce dependence among the
selected observations. The impact of such dependence on the empirical
cumulative distribution function (c.d.f.) is studied. An
asymptotic framework and weak conditions on the informative selection
mechanism are developed under which the (unweighted)
empirical c.d.f. converges uniformly, in $L_2$ and almost surely, to a
weighted version of the superpopulation c.d.f.
This yields an analogue of the Glivenko--Cantelli theorem. A series of
examples, motivated by real problems in
surveys and other observational studies, shows that the conditions are
verifiable for specified designs.
\end{abstract}

%

\end{frontmatter}
%

\section{Introduction}

Consider informative selection of a sample from a finite population,
with responses $Y$ realized as
independent and identically distributed (i.i.d.) random variables with
probability density function (p.d.f.) $f$,
referred to as the superpopulation model.
(Regression problems, in which observations are conditionally
independent given covariates, are also of interest, but the
following discussion readily generalizes to that setting and we
restrict attention to the i.i.d. case for simplicity of exposition.)
In non-informative selection (e.g., Cassel \textit{et al.} \cite
{CasselSarndalWretman1977}, Section
1.4, or S{\"a}rndal \textit{et al.}~\cite{SarndalSwenssonWretman92}, Remark
2.4.4),
the probability of drawing the sample does not depend explicitly on the
responses~$Y$.
We consider informative selection in the sense that the sample
responses, \textit{given that they were selected}, are not i.i.d. $f$.
A specification of informative selection that includes the i.i.d. case
described here is given in Pfeffermann and Sverchkov \cite
{PfeffermannSverchkov2009}, Remark 1.2.
We study the implications of this informative selection for estimation
of the superpopulation model.

In general, the informative selection mechanism may induce dependence
among the selected observations.
Nevertheless, a large body of current methodological literature treats
the observations as if they were independently distributed
according to the \textit{sample p.d.f.}, defined as the conditional
distribution of the random variable $Y$,
given that it was selected. Under informative selection, the sample p.d.f.
differs from $f$.
In particular, Pfeffermann \textit{et al.} \cite
{PfeffermannKriegerRinott98} (see some
motivating work in Skinner~\cite{Skinner1994}) have developed
a \textit{sample likelihood} approach to estimation and inference for the
superpopulation model, which maximizes the criterion function formed by taking
the product of the sample p.d.f.'s, as if the responses were i.i.d. This
methodology has been extended in a number of directions Eideh and Nathan
\cite{EidehNathan2007,EidehNathan2007c,EidehNathan2009},
Pfeffermann \textit{et al.}~\cite{PfeffermanMouraDaSilvaLuisdonascimento2006},
Pfeffermann and Sverchkov \cite
{PfeffermannSverchkov1999,PfeffermannSverchkov2003,PfeffermannSverchkov2007}.
An extensive review
of these and other approaches to inference under
informative selection is given by Pfeffermann and Sverchkov \cite
{PfeffermannSverchkov2009}.

Under a strong set of assumptions (in particular, sample size remains
fixed as population size goes to infinity), Pfeffermann \textit{et al.}
\cite{PfeffermannKriegerRinott98} have established the pointwise convergence
of the joint distribution of the responses to the product of the sample
p.d.f.'s. This is taken as partial justification of
the sample likelihood approach. Alternatively, the full likelihood of
the data (selection indicators for the finite population and response
variables and inclusion probabilities for the sample) can be maximized
(Breckling \textit{et~al.}~\cite{Brecklingetal1994}, Chambers \textit
{et al.}~\cite{MR1616053}), or the \textit{pseudo-likelihood} can
be obtained by plugging in Horvitz--Thompson estimators for unknown
quantities in the log-likelihood for the entire finite population
(e.g., Binder~\cite{Binder1983}, Chambers and Skinner \cite
{ChambersSkinner2003}, Kish and Frankel~\cite{KishFrankel1974},
Section~2.4).
Obviously, each of these likelihood-based approaches requires a model
specification.

Rather than starting at the point of likelihood-based inferential
methods for the superpopulation model,
we take a step back and consider the problem of identifying a suitable
model using observed data.
In an ordinary inference problem with i.i.d. observations,
we often begin not by constructing a likelihood and conducting
inference, but by using basic sample statistics to help identify a
suitable model.
In particular, under i.i.d. sampling the empirical cumulative
distribution function (c.d.f.) converges uniformly almost surely to the
population c.d.f.,
by the Glivenko--Cantelli theorem (e.g., van~der Vaart \cite
{vanderVaart1998}, Theorem~19.1). What is the behavior of the
empirical c.d.f. under informative selection
from a finite population? In this paper, we develop an asymptotic
framework and weak conditions on the informative selection mechanism
under which the (unweighted) empirical c.d.f. converges uniformly, in
$L_2$ and almost surely, to a weighted version of the superpopulation
c.d.f. The corresponding quantiles also converge uniformly on compact
sets. Our almost sure results rely on an embedding argument.
Importantly, our construction preserves the original response vector
for the finite population, not some independent replicate.

The conditions we propose are verifiable for specified designs, and
involve computing conditional versions of first and second-order
inclusion probabilities.
Motivated by real problems in surveys\vadjust{\goodbreak} and other observational studies,
we give examples of where these conditions hold and where they fail.
Where the conditions hold, the convergence results we obtain may be
useful in making inference about the superpopulation model.
For example, the results may be used in identifying a suitable
parametric family for the weighted c.d.f.,
from which a selection mechanism and a superpopulation p.d.f. may be
postulated using results in Pfeffermann \textit{et al.} \cite
{PfeffermannKriegerRinott98}.\vspace*{-2pt}

\section{Results}\vspace*{-1pt}
\subsection{Asymptotic framework and assumptions}\vspace*{-1pt}

In what follows, all random variables are defined on a common
probability space $(\Omega,\mathscr{A},P)$.
Let $\mathscr{B}(\mathbb{R})$ denote the $\sigma$-field of Borel sets.
Assume that for ${k \in \mathbb{N}}$, $Y_k \dvtx (\Omega,\mathscr
{A},\break P) \to(\mathbb{R},\mathscr{B}(\mathbb{R}))$
are i.i.d. real random variables with a density $f$ with respect to
$\lambda$, the Lebesgue measure. Consider
$\{N_\gamma\}_{\gamma\in\mathbb{N}}$, an increasing
sequence of positive integers representing a sequence of population
sizes, with $\lim_{\gamma\to\infty} N_\gamma=\infty$.

We consider a sequence of finite populations and samples. The $\gamma$th finite population is
the set of elements indexed by $U_\gamma=(1,\ldots,N_\gamma)$. In
the sampling literature (e.g., S{\"a}rndal \textit{et al.} \cite
{SarndalSwenssonWretman92}),
$U_\gamma$ is often an unordered set, but it is convenient for us to
order it and to write, for example, $\sum_{k\in
U_\gamma}=\sum_{k=1}^{N_\gamma}$. The vector of responses for the
population is $\mathcal{Y}_\gamma =(Y_k)_{k\in U_\gamma}$ and the
sample is
indexed by the random vector $\mathcal{I}_\gamma =(I_{\gamma
k})_{k\in U_\gamma}$,
where the $k$th coordinate $I_{\gamma k}$ indicates the number of
times element $k$ is selected: 0 or 1 under without-replacement
sampling, or a non-negative integer under with-replacement sampling.
Define the distribution of $\mathcal{I}_\gamma $ conditional
on~$\mathcal{Y}_\gamma $:
\[
g_\gamma(i_1\ldots,i_{N_\gamma},y_1,\ldots, y_{N_\gamma})=\mathrm
{P}\bigl(\mathcal{I}_\gamma
=(i_1,\ldots,i_{N_\gamma})\vert\mathcal{Y}_\gamma =(y_1,\ldots
,y_{N\gamma})\bigr).
\]

We assume that the index of the element $k$ of the population plays no
role in the way elements are selected. Specifically, let $\sigma$
denote a permutation of a vector of length $N_\gamma$. Then, for all
$\gamma\in\mathbb{N}$,
$(\mathcal{I}_\gamma \vert\mathcal{Y}_\gamma )$ and $(\sigma\cdot
\mathcal{I}_\gamma \vert\sigma\cdot\mathcal{Y}_\gamma )$ are
identically distributed, or equivalently
%
%
\begin{equation}
g_\gamma(i_1,\ldots, i_{N_\gamma},y_1,\ldots, y_{N_\gamma})
=g_\gamma\bigl(\sigma\cdot(i_1,\ldots, i_{N_\gamma}),\sigma\cdot
(y_1,\ldots
, y_{N_\gamma})\bigr).
\label{eq:exchangeability}
\end{equation}
We refer to (\ref{eq:exchangeability}) as the exchangeability
assumption. It corresponds to the condition of weakly exchangeable
arrays (Eagleson and Weber~\cite{EaglesonWeber1978}) applied to
$(I_{\gamma k},Y_k
)_{\gamma\in\mathbb{N},k\in U_\gamma}$.\vspace*{-1pt}

\begin{definition}\label{D1}
For $\gamma\in\mathbb{N}$, the empirical c.d.f. is the random process
$F_\gamma\dvtx \mathbb{R}\to[0,1]$ via
\[
F_\gamma(\alpha)=\frac{\sum_{k\in U_\gamma}\mathbh{1}_{(-\infty
,\alpha
]}(Y_k ) I_{\gamma k}}
{\mathbh{1}_{\mathcal{I}_\gamma =0}+\sum_{k\in U_\gamma}I_{\gamma k}}.\vspace*{-1pt}
\]
\end{definition}

\begin{definition}\label{D2}
Given $\gamma$, let $k,\ell\in U_\gamma$ with $k\ne\ell$. Assume
exchangeability as in (\ref{eq:exchangeability}) and~let
\begin{eqnarray*}
m_\gamma(y)&=& \mathrm{E}[I_{\gamma k}\vert Y_{k}=y],\\
v_\gamma(y)&=&\operatorname{Var}(I_{\gamma k}\vert Y_{k}=y),\vadjust{\goodbreak}\\
m'_\gamma(y_1,y_2)&=&\mathrm{E}[I_{\gamma k}\vert Y_{k}=y_1,Y_{\ell
}=y_2],\\
c_\gamma(y_1,y_2)&=&\operatorname{Cov}(I_{\gamma k},I_{\gamma\ell
}\vert
Y_{k}=y_1,Y_{\ell}=y_2).
\end{eqnarray*}
(These definitions do not depend on the choice of
$k, \ell$ under the exchangeability assumption).
\end{definition}

The following conditions on $m_\gamma$ are used in defining the limit c.d.f.:
\setcounter{assumption}{-1}
\setcounter{equation}{-1}
\begin{assumption} There exist $M\dvtx\mathbb{R}\to\mathbb{R}^+$ and
$m\dvtx\mathbb{R}\to\mathbb{R}^+$, both $\lambda$-measurable, such that
\begin{subequation}\label{A0}
%
%
\begin{equation}\label{A0.1}
\cases{
\displaystyle\forall\gamma\in\mathbb{N},\qquad m_\gamma
<M,\vspace*{2pt}\cr
\displaystyle\int Mf\, \mathrm{d}\lambda<\infty,}
\end{equation}\vspace*{-6pt}
%
%
\begin{equation}\label{A0.2}
\cases{
m_\gamma\to m \qquad\mbox{pointwise as $\gamma\to\infty$},\vspace
*{2pt}\cr
\displaystyle\int mf \,\mathrm{d}\lambda>0.}
\end{equation}
\end{subequation}
\end{assumption}

\begin{definition}
\label{D3}
Under \textup{A\ref{A0}}, the limit c.d.f. $F_s\dvtx\mathbb{R}\to[0,1] $ is
\[
F_s(\alpha)=\frac{\int\mathbh{1}_{(-\infty,\alpha]} mf \,\mathrm
{d}\lambda
}{\int mf \,\mathrm{d}\lambda}.
\]

\end{definition}

\begin{remark*}[(relation to sample p.d.f.)] Because of informative
selection, the empirical c.d.f. does not converge to the superpopulation
c.d.f. Under some conditions to be specified below, it converges to
$F_s$, a weighted integral of the superpopulation p.d.f. To see this,
consider the case of without-replacement sampling and a single element,
$k$. The sample p.d.f. defined in Krieger and Pfeffermann \cite
{KriegerPfeffermann1992} is the
conditional density of $Y_k $ given $I_{\gamma k}=1$. By Bayes' rule,
\begin{eqnarray*}
f_{s\gamma}(y)&=&f(y\vert I_{\gamma k}=1)=\frac{\mathrm{P}(I_{\gamma
k}=1\vert Y_k =y)f(y)}{\int\mathrm{P}(I_{\gamma k}=1\vert Y_k =y)f
\,\mathrm{d}\lambda}\\
&=&\frac{m_{\gamma}(y)}{\int m_\gamma f \,\mathrm{d}\lambda
}f(y)=w_\gamma(y)f(y).
\end{eqnarray*}
Define $w=\lim_{\gamma\to\infty}w_\gamma$ and consider $\alpha\in
\mathbb
{R}$. Then
\[
\lim_{\gamma\to\infty}\int\mathbh{1}_{(-\infty,\alpha]}
f_{s\gamma}
\,\mathrm{d}\lambda=
\lim_{\gamma\to\infty}\int\mathbh{1}_{(-\infty,\alpha]}
w_\gamma f
\,\mathrm{d}\lambda=
\int\mathbh{1}_{(-\infty,\alpha]} w f \,\mathrm{d}\lambda
=F_s(\alpha).
\]
Thus, if observations were i.i.d. from the sample p.d.f., $F_s$ would be
the natural limiting c.d.f. A related argument can be used to show that
the same weighted c.d.f. is obtained under with-replacement sampling
and a
fixed number of draws, when considering the distribution of any
observation in the sample.

Because informative selection from a finite population may induce
dependence among the selected observations, observations are not
i.i.d., and we next specify asymptotic weak dependence conditions among
$\mathcal{I}_\gamma $ coordinates.
\end{remark*}

For a sequence $\{b_\gamma\}$, let $\mathrm{o}_\gamma(b_\gamma)$
denote $\lim
_{\gamma\to\infty}\mathrm{o}_\gamma(b_\gamma)b_\gamma^{-1}=0$.
In the next two
assumptions, we define sufficient conditions for uniform $L_2$
convergence and uniform a.s. convergence of the empirical c.d.f.

\begin{assumption}[(Uniform $L_2$ convergence conditions)]
%
\begin{subequation}\label{A1}
%
%
\begin{eqnarray}
\int c_\gamma(y_1,y_2)f(y_1)f(y_2)\,\mathrm{d}y_1 \,\mathrm{d}y_2
&=&\mathrm{o}_\gamma(1), \label{A1.1}
\\
\int\bigl(m_\gamma'(y_1,y_2)m'_\gamma(y_2,y_1)-m_\gamma(y_1)m_\gamma
(y_2)\bigr)f(y_1)f(y_2)\,\mathrm{d}y_1 \,\mathrm{d}y_2 &=&\mathrm
{o}_\gamma(1), \label{A1.2}
\\
\int(v_\gamma+m_\gamma^2)f \,\mathrm{d}\lambda&=&\mathrm
{o}_\gamma(N_\gamma),\label{A1.3}
\\
\mathrm{P}\bigl(\mathcal{I}_\gamma =(0,\ldots, 0)\bigr)&=&\mathrm
{o}_\gamma(1).\label{A1.5}
\end{eqnarray}
\end{subequation}
\end{assumption}

\begin{assumption}[(Uniform almost sure convergence conditions)]
Let $y \in\mathbb{R}^\mathbb{N}$ satisfy
\[
\sup_{\alpha'\in\mathbb{R}}\biggl|\frac{\sum_{k\in U_\gamma
}\mathbh{1}_{(-\infty,\alpha']}(y_k)}{N_\gamma}-
\int\mathbh{1}_{(-\infty,\alpha']} f \,\mathrm{d}\lambda\biggr
|=\mathrm{o}_\gamma(1).
\]
Then for all $\alpha\in\mathbb{R}$,
\begin{subequation}\label{A2}
%
%
\begin{eqnarray}
\operatorname{Var}\biggl(\sum_{k\in U_\gamma}\mathbh{1}_{(-\infty
,\alpha
]}(y_k)I_{\gamma k}\big| \mathcal{Y}_\gamma =(y_1,\ldots,
y_{N_\gamma})\biggr)&=&\mathrm{o}_\gamma
(N_\gamma^2),\label{A2.1}
\\
\sum_{k\in U_\gamma} \mathbh{1}_{(-\infty,\alpha]}(y_k)
\bigl(\mathrm{E}[I_{\gamma k}\vert\mathcal{Y}_\gamma =(y_1,\ldots
, y_{N_\gamma})]-m_{\gamma
}(y_k)\bigr)
&=&\mathrm{o}_\gamma(N_\gamma), \label{A2.2}
\\
g_\gamma((0,\ldots,0),y)&=&\mathrm{o}_\gamma(1). \label{A2.3}
\end{eqnarray}
\end{subequation}
%
\end{assumption}

\subsubsection*{Properties of sampling without replacement}

In the case of sampling without replacement, $\mathcal{I}_\gamma
\dvtx \Omega\to
\{0,1\}^{N\gamma}$, \textup{A\ref{A0}} and \textup{A\ref{A1}} can be replaced by a
simpler set of sufficient conditions for uniform $L_2$ convergence.

\begin{assumption}[(Uniform $L_2$ convergence conditions under
sampling without replacement)] \label{A3}
\begin{subequation}
%
%
\begin{equation}\label{A3.2}
\exists m\dvtx\mathbb{R}\to\mathbb{R}^+ \ \lambda\mbox{-measurable
s.t.}\qquad
\cases{
m_\gamma\to m \qquad\mbox{pointwise as
$\gamma\to\infty$},\vspace*{2pt}\cr
\displaystyle\int mf \,\mathrm{d}\lambda>0,}
\end{equation}
%
%
\begin{eqnarray}
\forall y_1,y_2,\qquad c_\gamma(y_1,y_2)&=&\mathrm{o}_\gamma(1),
\label{A3.3}
\\
\forall y_1,y_2,\qquad m_\gamma'(y_1,y_2)-m_\gamma(y_2)&=&\mathrm
{o}_\gamma(1), \label{A3.4}
\\
\mathrm{P}\bigl(\mathcal{I}_\gamma =(0,\ldots, 0)\bigr)&=&\mathrm
{o}_\gamma(1).\label{A3.5}
\end{eqnarray}
These conditions imply \textup{A\ref{A0}} and \textup{A\ref{A1}}.
\end{subequation}
\end{assumption}

\begin{pf}
Since $I_{\gamma k}\in\{0,1\}$, (\ref{A0.1}) and (\ref{A1.3}) always
hold. By applying the Lebesgue dominated convergence theorem,
we obtain that (\ref{A1.1}) is verified when $\forall y_1,y_2,
c_\gamma
(y_1,y_2)=\mathrm{o}_\gamma(1)$
and (\ref{A1.2}) is verified when $\forall y_1,y_2,\  m_\gamma
'(y_1,y_2)-m_\gamma(y_2)=\mathrm{o}_\gamma(1)$.
\end{pf}

An important special case of sampling without replacement is
non-informative selection, with~$\mathcal{I}_\gamma $ independent of
$\mathcal{Y}_\gamma $ for all
$\gamma\in\mathbb{N}$. In this case, the sample obtained is an i.i.d.
sample of size $n_\gamma=\sum_{k\in U_\gamma}I_{\gamma k}$
(Fuller~\cite{Fuller2009}, Theorem~1.3.1), and the classic
Glivenko--Cantelli theorem can be applied as soon as
$n_\gamma{\mathop{\rightarrow}^{\mathrm{a.s.}}}\infty$ as $\gamma
\to
\infty$. The assumptions of Theorem~\ref{T1} and Theorem~\ref{T2} will
then just ensure that the expectation of the sample size will grow to infinity,
and that its variations are small enough to avoid very small samples.
We can thus replace A\ref{A0}--A\ref{A2} by a simpler set of sufficient
conditions.

\begin{assumption}[(Uniform $L_2$ and a.s. convergence
conditions under independent sampling without replacement)]\label{A4}
%
%
\begin{equation} \cases{
{N_\gamma}^{-1}\mathrm{E}[n_\gamma]\to m\ne0\qquad\mbox{as
$\gamma\to\infty$},\vspace*{2pt}\cr
\operatorname{Var}(n_\gamma)=\mathrm{o}_\gamma(N_\gamma^2).
}
\end{equation}
These conditions imply \textup{A\ref{A0}}--\textup{A\ref{A2}}.
\end{assumption}

\begin{pf}
We first show that \textup{A\ref{A4}} implies \textup{A\ref{A3}}.
Because $\mathcal{I}_\gamma $ and $\mathcal{Y}_\gamma $ are
independent, the exchangeability
assumption implies
$m_\gamma(y) = \mathrm{E}[I_{\gamma1}] = N_\gamma^{-1}\mathrm
{E}[n_\gamma]$
and ${N_\gamma}^{-1}\mathrm{E}[n_\gamma]\to m$ by~A\ref{A4}, so
(\ref{A3.2})
holds. Exchangeability also implies
\[
\mathrm{E}[I_{\gamma1}I_{\gamma2}] = \frac{\sum_{k,\ell\in
U_\gamma\dvtx k\ne
\ell}\mathrm{E}[I_{\gamma k}I_{\gamma\ell}]}{N_\gamma(N_\gamma-1)}
= \mathrm{E}\biggl[\frac{\sum_{k,\ell\in U_\gamma\dvtx k\ne
\ell}I_{\gamma k}I_{\gamma
\ell}}{N_\gamma(N_\gamma-1)}\biggr]
= \mathrm{E}\biggl[\frac{n_\gamma(n_\gamma- 1)}{N_\gamma
(N_\gamma- 1)} \biggr]
\]
so
%
%
\begin{equation}
c_\gamma(y_1,y_2)=\operatorname{Cov}(I_{\gamma1},I_{\gamma2})
=\mathrm{E}\biggl[\frac{n_\gamma(n_\gamma-N_\gamma)}{N_\gamma
^2(N_\gamma-1)}\biggr]
+\operatorname{Var}\biggl(\frac{n_\gamma}{N_\gamma}\biggr
)=\mathrm{o}_\gamma(1)\label{eq:cov}
\end{equation}
by \textup{A\ref{A4}}, so (\ref{A3.3}) is obtained, and (\ref{A3.4}) holds by
independence. Finally,
%
%
\begin{eqnarray}\label{eq:Pn0}
\mathrm{P}(n_\gamma=0)&=&\mathrm{P}(n_\gamma<1)=\mathrm
{P}(n_\gamma-\mathrm{E}[n_\gamma]<1-\mathrm{E}[n_\gamma
])
\nonumber
\\[-8pt]
\\[-8pt]
\nonumber
&\leq& \mathrm{P}(|n_\gamma-\mathrm{E}[n_\gamma]|>\mathrm
{E}[n_\gamma]-1)
\leq\frac{\operatorname{Var}(n_\gamma)}{(\mathrm{E}[n_\gamma]-1)^2}
= \mathrm{o}_\gamma(1),
\end{eqnarray}
establishing (\ref{A3.5}).

We next show that (\ref{A4}) implies \textup{A\ref{A2}}. For all $\alpha\in
\mathbb{R}$,
\begin{eqnarray*}
&&\operatorname{Var}\biggl(\sum_{k\in U_\gamma} \mathbh
{1}_{(-\infty,\alpha
]}(Y_k)I_{\gamma k}\big| \mathcal{Y}_\gamma =(y_1,\ldots,
y_{N_\gamma})\biggr)\\
&&\quad=\sum_{k\in U_\gamma} \mathbh{1}_{(-\infty,\alpha
]}(y_k)\operatorname{Var}(I_{\gamma
k})\\
&&\qquad{}+\sum_{k,\ell\in U_\gamma\dvtx k\ne\ell} \mathbh
{1}_{(-\infty,\alpha
]}(y_k)\mathbh{1}_{(-\infty,\alpha]}(y_\ell)\operatorname
{Cov}(I_{\gamma
k},I_{\gamma\ell})\\
&&\quad\le N_\gamma+N_\gamma(N_\gamma-1)\mathrm{o}_\gamma
(1)=\mathrm{o}_\gamma(N_\gamma
^2)
\end{eqnarray*}
by equation (\ref{eq:cov}), so (\ref{A2.1}) holds. By independence,
\[
\mathrm{E}[I_{\gamma k}\vert\mathcal{Y}_\gamma =(y_1,\ldots
,y_{N_\gamma})]=\mathrm{E}[I_{\gamma k}\vert
Y_k=y_k]=m_\gamma(y_k),
\]
so (\ref{A2.2}) holds. Finally,
\[
g_\gamma((0,\ldots,0),y)=\mathrm{P}(n_\gamma=0)=\mathrm{o}_\gamma(1)
\]
by independence and (\ref{eq:Pn0}), so (\ref{A2.3}) holds.
\end{pf}

\begin{remark*} In conventional finite population asymptotics
(Breidt and Opsomer~\cite{bre00,bre08b}, Isaki and Fuller \cite
{IsakiFuller82}, Robinson and S\"{a}rndal~\cite{rob83}),
conditions on design
covariances Cov$(I_{\gamma k},I_{\gamma\ell})$ are imposed to
guarantee that the Horvitz--Thompson estimator $\sum_{k\in U_\gamma
}y_kI_{\gamma k}(\mathrm{E}[I_{\gamma k}])^{-1}$ is consistent. Typically,
these conditions imply that the variance of the Horvitz--Thompson
estimator is $\mathrm{O}_\gamma(N_\gamma^2/(N_\gamma\pi_{*\gamma
}))$, where
$N_\gamma\pi_{*\gamma}\to\infty$ is a sequence of lower bounds on the
expected sample size, $\mathrm{E}[n_\gamma]$. These same conditions
can be used
to show that $\operatorname{Var}(n_\gamma)= \mathrm{O}_\gamma
(N_\gamma^2/(N_\gamma\pi_{*\gamma
}))=\mathrm{o}_\gamma(N_\gamma^2)$, agreeing with \textup{A\ref{A4}}.
\end{remark*}
%
\subsection{Uniform convergence of the empirical c.d.f.}

In this section, we state the main results of the paper: uniform $L_2$
convergence of the empirical c.d.f. and uniform almost sure convergence of
the empirical c.d.f. Important corollaries yield uniform convergence of
sample quantiles on compact sets. Proofs are given in the \hyperref
[appa]{Appendix}.


\subsubsection{Uniform $L_2$ convergence of the empirical c.d.f.}

\begin{theorem}\label{T1} Under \textup{A\ref{A0}} and \textup{A\ref{A1}}, the empirical
c.d.f. converges uniformly in $L_2$ in the sense that
\[
\sup_{\alpha\in\mathbb{R}}|F_\gamma(\alpha)-F_s(\alpha) |=\|
F_\gamma-F_s\|_\infty \mathop{\rightarrow}\limits_{\gamma\to\infty}^{{L_2}}0.
\]
\end{theorem}

\begin{definition} The limit quantiles $\xi_s\dvtx(0,1)\to\mathbb
{R}$ are
given by
\[
\xi_s(p)=\inf\{y\in\mathbb{R}\dvtx F_s(y)\geq p\}
\]
and the empirical quantiles $\xi_\gamma\dvtx(0,1)\to\mathbb{R}$
are given by
\[
\xi_\gamma(p)=\inf\{y\in\mathbb{R}\dvt F_\gamma(y)\geq p\}.
\]
\end{definition}

With this definition, we have the following corollary.
\begin{corollary}\label{C1}
Suppose that $F_s$ is continuous on $\mathbb{R}$ and
$0<F_s(y_1)=F_s(y_2)<1\Rightarrow y_1=y_2$. Then, under
\textup{A\ref{A0}} and \textup{A\ref{A1}}, the empirical quantiles converge uniformly in
probability to the limit quantiles,
\[
\sup_{p\in K}| \xi_\gamma(p)-\xi_s(p)|\mathop{\rightarrow}\limits_{\gamma\to
\infty}^{{P}}0
\]
for all $K$ a compact subset of $(0,1)$.
Under the further hypothesis that $f$ has compact support, the
convergence is uniform in $L_2$:
\[
\sup_{p\in K}| \xi_\gamma(p)-\xi_s(p)|\mathop{\rightarrow}\limits_{\gamma\to
\infty}^{{L_2}}0.
\]

\end{corollary}


\subsubsection{Uniform almost sure convergence of the empirical c.d.f.}

The Glivenko--Cantelli theorem gives uniform almost sure convergence of
the empirical c.d.f. under i.i.d. sampling. We now consider uniform almost
sure convergence under dependent sampling satisfying the second-order
conditions of A\ref{A2}.

Asymptotic arguments in survey sampling consist first in embedding
a specific sample scheme in a sequence of sample schemes.
In the proof of the following representation theorem, we link the
elements of the sequence of sample schemes in a way
that ensures uniform almost sure convergence of the empirical
c.d.f. We stress that in our result the vector of responses for the
population remains the original $\mathcal{Y}_\gamma =(Y_k)_{k\in
U_\gamma}$, and
not another set of identically distributed random variables.

\begin{theorem}\label{T2}
Under \textup{A\ref{A0}} and \textup{A\ref{A2}}, there exist sequences of random variables
$(I'_{\gamma k})_{\gamma\in\mathbb{N},k\in U_\gamma}$,
$(Y'_k)_{k\in\mathbb{N}}$ defined on
the probability space $(\Omega\times[0,1],\mathscr{A}\otimes
\mathscr{B}_{[0,1]},\mathrm
{P}'=\mathrm{P}\otimes\lambda_{[0,1]})$
such that
\begin{itemize}
\item$\|F'_\gamma-F_s \|_\infty$ converges $\mathrm{P}'$-a.s. to $0$
\item$\forall\gamma\in\mathbb{N}, (\mathcal{I}_\gamma ',\mathcal
{Y}_\gamma ')$ and $(\mathcal{I}_\gamma ,\mathcal{Y}_\gamma )$ have
the same law
\item$\forall\gamma\in\mathbb{N}, \omega\in\Omega, x\in[0,1],\
\mathcal{Y}_\gamma
'(\omega,x)=\mathcal{Y}_\gamma (\omega)$,
\end{itemize}
where $\mathscr{B}_{[0,1]}$ is the $\sigma$-field of Borel sets,
$\lambda_{[0,1]}$ is the Lebesgue measure on $[0,1]$,
$\mathcal{I}_\gamma '=(I'_{\gamma1},\ldots, I'_{\gamma N_\gamma
})$, $\mathcal{Y}_\gamma
'=(Y'_1,\ldots, Y'_{N_\gamma})$
and $F'_\gamma\dvtx \mathbb{R}\to[0,1]$ via
%
%
\begin{equation}
F'_\gamma(\alpha)= \frac{\sum_{k\in U_\gamma}\mathbh{1}_{(-\infty
,\alpha]}(Y'_{\gamma k})I'_{\gamma k}}
{\sum_{k\in U_\gamma}I'_{\gamma k}+\mathbh{1}_{I'_\gamma=0}}.
\label{eq:Fprime}
\end{equation}
\end{theorem}

\begin{corollary}\label{C2} Suppose that $F_s$ is continuous and
$0<F_s(y_1)=F_s(y_2)<1\Rightarrow y_1=y_2$. If \textup{A\ref{A0}} and \textup{A\ref{A2}}
hold, then for $(I'_{\gamma k})_{\gamma\in\mathbb{N},k\in U_\gamma
}$ and
$(Y'_k)_{k\in\mathbb{N}}$ that satisfy the conditions of Theorem \ref
{T2}, the empirical quantiles
\[
\xi'_\gamma(p)=\inf\{y\in\mathbb{R}\dvt F'_\gamma(y)\geq p\}
\]
converge uniformly almost surely,
\[
\sup_{p\in K}| \xi_\gamma(p)-\xi_s(p)|\mathop{\rightarrow}\limits_{\gamma\to
\infty}^{{\mathrm{a.s.}}}0
\]
for all $K$ a compact subset of $(0,1)$.

\end{corollary}
%
\section{Examples}

We now consider a series of examples of selection mechanisms,
motivated by real problems in surveys and other observational studies.
We give examples where conditions A\ref{A0}, A\ref{A1}, A\ref{A2} hold
and where they fail.

\subsection{Non-informative selection without replacement}

\begin{itemize}
\item For any sequence of \textit{fixed-size without-replacement
designs} with $\mathcal{I}_\gamma $ independent of $\mathcal
{Y}_\gamma $ (e.g., simple random sampling,
stratified sampling
with stratification variables independent of $\mathcal{Y}_\gamma $,
rejective sampling (H{\'a}jek~\cite{Hajek81}) with inclusion probabilities
independent of $\mathcal{Y}_\gamma $, etc.), the condition A\ref{A4}
holds provided
that $n_\gamma N_\gamma^{-1}$
converges to a strictly positive sampling rate.

\item For a sequence of \textit{Bernoulli samples} with parameter
$p\in(0,1)$, the $\{I_{\gamma k}\}$ are i.i.d. $\operatorname{Bernoulli}(p)$ random
variables, so $\operatorname{Var}(n_\gamma)=N_\gamma p(1-p)$
and condition A\ref{A4} holds.

\item\textit{Poisson sampling} corresponds to a design in which, given
a random vector
$(\Pi_{\gamma1},\ldots, \Pi_{\gamma N_\gamma})\dvtx \Omega\to
(0,1]^{N_\gamma}$,
the $\{I_{\gamma k}\}$ are a sequence of independent Bernoulli($\Pi
_{\gamma k}$) random variables (Poisson~\cite{Poisson1837}).
In this case, the variance of $n_\gamma$ is given by
\[
\operatorname{Var}(n_\gamma) =\sum_{k\in U_\gamma} \mathrm{E}[\Pi
_{\gamma k}(1-\Pi_{\gamma
k})]+\operatorname{Var}\biggl(\sum_{k\in U_\gamma} \Pi_{\gamma k}\biggr).
\]
Note that the first term in this expression is always $\mathrm
{o}_\gamma(N_\gamma
^2)$, so it suffices to consider the second.

\begin{itemize}[--]
\item[--] In the case where the vector $[\Pi_{\gamma k}]_{k\in U_\gamma}$
is just a random permutation of a non-random vector $[\pi_{\gamma
k}]_{k\in U_\gamma}$, then $\operatorname{Var}(\sum_{k\in U_\gamma
}\Pi_{\gamma k})=\operatorname{Var}(\sum_{k\in U_\gamma}\pi
_{\gamma k})=0$
and A\ref{A4} is satisfied when
$ {N_\gamma}^{-1} \sum_{k\in U_\gamma}\pi_{\gamma k}$ converges to a
non-zero constant.

\item[--] Suppose that $Z_{\gamma}$ is a random positive real vector of
size $N_\gamma$, and suppose that the law of $(Z_\gamma,\mathcal
{Y}_\gamma )$ is
invariant under any permutation of the coordinates. For
$n^*_\gamma$
fixed, consider the design in which $\Pi_{\gamma k}=n^*_\gamma
{Z_{\gamma k}}\{\sum_{k\in U_\gamma} Z_{\gamma k}\}^{-1}$. Then
\[
\operatorname{Var}\biggl(\sum_{k\in U_\gamma}\Pi_{\gamma
k}\biggr)=\operatorname{Var}(n_\gamma^*)=0
\]
and A\ref{A4} is satisfied when $Z_\gamma$ and $\mathcal{Y}_\gamma
$ are independent and
${N_\gamma}^{-1}{n^*_\gamma}$ converges to a non-zero constant.

\item Let $a_\gamma,b_\gamma\in(0,1]$ with $a_\gamma\ne b_\gamma$. If
\[
(\Pi_{\gamma1},\ldots,\Pi_{\gamma N_\gamma})\equiv
\cases{
(a_\gamma,\ldots, a_\gamma), & \quad$\mbox{with probability
$1/2$},$\vspace*{2pt}\cr
(b_\gamma,\ldots, b_\gamma), & $\quad\mbox{with probability $1/2$},$}
\]
then
\[
\operatorname{Var}\biggl(\sum_{k\in U_\gamma}\Pi_{\gamma k}\biggr)=N_\gamma
^2\frac{(a_\gamma
-b_\gamma)^2}{4}\ne\mathrm{o}_\gamma(N_\gamma^2).
\]
Then A\ref{A4} is not verified and in fact if $N_\gamma a_\gamma
=\mathrm{o}_\gamma(1)$ we do not have uniform convergence of the
empirical c.d.f.
\end{itemize}
\end{itemize}

\subsection{Length-biased sampling}

Length-biased sampling, in which $\mathrm{P}(I_{\gamma k}=1\vert
Y_k=y_k)=m_\gamma(y_k)\propto y_k$, is pervasive in real surveys and
observational studies. Cox~\cite{cox69} gives a now-classic example of
sampling fibers in textile manufacture, in which $m_\gamma(y_k)\propto
y_k=$ fiber length. In surveys of wildlife abundance, ``visibility
bias'' means that larger individuals or groups are more noticeable
(e.g., Patil and Rao~\cite{PatilRao1978}), so $m_\gamma(y_k)\propto
y_k=$ size of
individual or group. ``On-site surveys'' are sometimes used to study
people engaged in some activity like shopping in a mall
(Nowell and Stanley~\cite{NowellStanley1991}) or fishing at the seashore
(Sullivan \textit{et~al.}~\cite{Sullivan2006});
the longer they spend doing the activity, the more likely the field
staff are to intercept and interview them, so $m_\gamma(y_k)\propto
y_k=$ activity time. In mark-recapture surveys of wildlife populations,
individuals that live longer are more likely to be recaptured, so
$m_\gamma(y_k)\propto y_k=$ lifetime (e.g., Leigh~\cite{Leigh1988}).
Similarly, in epidemiological studies of latent diseases, individuals
who become symptomatic seek treatment and drop out of eligibility for
sampling, while those with long latency periods are more likely to be
sampled: $m_\gamma(y_k)\propto y_k=$ latency period. Finally,
propensity to respond to a survey is often related to a variable of
interest; for example, higher response rates from higher-income individuals.

Suppose that $f$ has compact, positive support: $\int\mathbh
{1}_{[\epsilon,M]}f \,\mathrm{d}\lambda=1$
for some $0<\epsilon<M<\infty$. For the $\gamma$th finite population,
consider Poisson sampling with inclusion probability proportional to $Y$,
in the sense that $\{I_{\gamma k}\}_{k\in U_\gamma}$ are independent
binary random variables, with
\[
\mathrm{P}(I_{\gamma k}=1\vert Y_k=y_k)=1-\mathrm{P}(I_{\gamma
k}=0\vert Y_k=y_k)=m_\gamma
(y_k)\propto y_k.
\]
Let $\tau_\gamma=y_k^{-1} \mathrm{P}(I_{\gamma k}=1\vert Y_k=y_k)$
be the common
proportionality constant (independent of $k$), and
assume that $\tau_\gamma\to\tau\in(0,M^{-1}]$ as $\gamma\to
\infty$. Then
\begin{eqnarray*}
m_\gamma(y)&=&\tau_\gamma y\quad \to \quad\tau y =m(y),\\
c_\gamma(y_k,y_\ell)&=&0, \qquad
m'_\gamma(y_k,y_\ell)-m_\gamma(y_k)=0,\\
\mathrm{P}\bigl(\mathcal{I}_\gamma =(0,\ldots,0)\bigr)&=&\mathrm{E}\biggl[\prod
_{k\in U_\gamma}(1-\tau_\gamma y_k)\biggr]\\
&\leq&(1-\tau_\gamma\epsilon)^{N_\gamma}=\exp\bigl(N_\gamma\ln
(1-\tau
_\gamma\epsilon)\bigr)=\mathrm{o}_\gamma(1),
\end{eqnarray*}
so that A\ref{A3} is verified. It then follows that the limiting
c.d.f. is
given by
%
%
\begin{equation}
F_s(\alpha)=\int\mathbh{1}_{(-\infty,\alpha]}\frac{y}{\mathrm{E}[Y_1]}f
\,\mathrm{d}\lambda. \label{eq:LBS}
\end{equation}

\subsection{Cluster sampling}
Let $F$ denote the superpopulation c.d.f.: $F(\tau)=\int\mathbh
{1}_{(-\infty,\tau]} f \,\mathrm{d}\lambda$.
Let $\tau\in\mathbb{R}$ be such that \mbox{$F(\tau)>0$}. Define
$i_{1\gamma
}=(\mathbh{1}_{(-\infty,\tau]}(Y_k))_{k\in U_\gamma}$ and
$i_{2\gamma
}=(\mathbh{1}_{(\tau,\infty)}(Y_k))_{k\in U_\gamma}$. The selection
mechanism is $\mathcal{I}_\gamma =i_{1\gamma}$ or $i_{2\gamma}$,
each with probability
$1/2$, so uniform convergence of the empirical c.d.f. is not possible.
Note that
\begin{eqnarray*}
\operatorname{Cov}(I_{\gamma k},I_{\gamma\ell}\vert
Y_{k}=y_1,Y_{\ell}=y_2)
&=&
\tfrac{1}{2}\mathbh{1}_{(-\infty,\tau]}(y_1)\mathbh{1}_{(-\infty
,\tau
]}(y_2)\\
&&{}+\tfrac{1}{2}\mathbh{1}_{(\tau,\infty)}(y_1)\mathbh{1}_{(\tau
,\infty)}(y_2)
-\tfrac{1}{4}
\end{eqnarray*}
so that
\[
\int c_\gamma(y_1,y_2)f(y_1)f(y_2) \,\mathrm{d}y_1 \,\mathrm{d}y_2=
\frac{1}{2}F^2(\tau)+\frac{1}{2}\bigl(1-F(\tau)\bigr)^2-\frac{1}{4}\ne\mathrm{o}_\gamma(1),
\]
and (\ref{A1.1}) fails to hold.
This example can be regarded as a ``worst-case'' cluster sample: the
sample consists of many elements but only one cluster, and the
population is made up of a small number of large clusters, none of
which is fully representative of the population.

\subsection{Cut-off sampling and take-all strata}

In cut-off sampling a part of the population is excluded from sampling,
so that $I_{\gamma k}=0$ with probability one for some subset of
$U_\gamma$. This may be due to physical limitations of the sampling
apparatus, like a net that lets small animals escape, or may be due to
a deliberate design decision. For example, a statistical agency may be
willing to accept the bias inherent in cutting off small $y$-values if
the $y$-distribution is highly skewed, as is often the case in
establishment surveys (e.g., S{\"a}rndal \textit{et al.}
\cite{SarndalSwenssonWretman92},
Section 14.4).

Consider cut-off sampling with $I_{\gamma k}=0$ for $\{k\in U_\gamma
\dvt
y_k\le\tau\}$, and simple random sampling without replacement of size
$\min\{n_\gamma, N_\gamma- \sum_{j\in U_\gamma} \mathbh
{1}_{(-\infty
,\tau]}(y_j ) \}$ from the remaining population, $\{j\in U_\gamma
\dvt
y_j>\tau\}$.

Define $Z_k=\mathbh{1}_{(-\infty,\tau]}(Y_k)$ with corresponding
realization $z_k=\mathbh{1}_{(-\infty,\tau]}(y_k)$.
Let $\rho_\gamma= N_\gamma^{-1}n_\gamma$ and assume that $\lim
_{\gamma
\to\infty} \rho_\gamma=\rho.$ We now verify A\ref{A3}.

Define $S_{\gamma A}=\sum_{j\in U_\gamma\dvtx j\notin A}Z_j$. By the weak
law of large numbers, $N_\gamma^{-1}S_{\gamma A}{\to}^{P}
F(\tau
)$ as $\gamma\to\infty$ for $A=\{k\}$ or $A=\{k,\ell\}$, and so for
those sets $A$ we have
\[
\lim_{\gamma\to\infty}
\mathrm{E}\biggl[\frac{\rho_\gamma-{N_\gamma}^{-1}S_{\gamma
A}}{1-N_\gamma
^{-1}S_{\gamma A}}
\mathbh{1}_{\{\rho_\gamma>{N_\gamma}^{-1}S_{\gamma A}\}}
\biggr]
=\frac{(\rho-F(\tau))\mathbh{1}_{\{\rho>F(\tau)\}}}{1-F(\tau)}
\]
by the uniform integrability of the integrand. With the same argument,
\begin{eqnarray*}
&& \lim_{\gamma\to\infty}\mathrm{E}\biggl[
\frac{(n_\gamma-S_{\gamma\{k,\ell\}})(n_\gamma
-1-S_{\gamma\{k,\ell\}})}
{(N_\gamma-S_{\gamma\{k,\ell\}})(N_\gamma-1-S_{\gamma\{
k,\ell\}})}
\mathbh{1}_{\{n_\gamma>S_{\gamma\{k,\ell\}}\}}\biggr]\\
&&\quad =\biggl(\frac{\rho-F(\tau)}{1-F(\tau)}\biggr)^2\mathbh{1}_{\{\rho
>F(\tau)\}}.
\end{eqnarray*}

Using conditional first and second-order inclusion probabilities under
simple random sampling, we have
\begin{eqnarray*}
m_\gamma(y_k)&=& z_k+(1-z_k)\mathrm{E}\biggl[\frac{n_\gamma-S_{\gamma\{k\}
}}{N_\gamma
-S_{\gamma\{k\}}}
\mathbh{1}_{\{n_\gamma>S_{\gamma\{k\}}\}}\biggr]\\
&\to& z_k+(1-z_k)\frac{(\rho-F(\tau))\mathbh{1}_{\{\rho>F(\tau)\}
}}{1-F(\tau)},\\
m'_\gamma(y_\ell,y_{k})&=&z_k+(1-z_\ell)(1-z_k)\mathrm{E}\biggl[\frac
{n_\gamma
-S_{\gamma\{k,\ell\}}}{N_\gamma-S_{\gamma\{k,\ell\}}}
\mathbh{1}_{\{n_\gamma>S_{\gamma\{k,\ell\}}\}}\biggr]\\
&&{} + z_\ell(1-z_k)\mathrm{E}\biggl[\frac{n_\gamma-1-S_{\gamma\{k,\ell\}
}}{N_\gamma
-1-S_{\gamma\{k,\ell\}}}
\mathbh{1}_{\{n_\gamma-1>S_{\gamma\{k,\ell\}}\}}\mathbh{1}_{\{
N_\gamma
-1>S_{\gamma\{k,\ell\}}\}}\biggr]\\
&\to& z_k+(1-z_k)\frac{(\rho-F(\tau))\mathbh{1}_{\{\rho>F(\tau)\}
}}{1-F(\tau)},\\
d_\gamma(y_k,y_\ell)&=&\mathrm{E}[I_{\gamma k}I_{\gamma\ell}\vert
Y_k=y_k,Y_\ell
=y_\ell]\\[-0.5pt]
&=&z_k z_\ell+\{z_k(1-z_\ell)+(1-z_k)z_\ell\}\mathrm{E}\biggl[\frac
{n_\gamma-1-S_{\gamma\{k,\ell\}}}{N_\gamma-1-S_{\gamma\{k,\ell\}}}
\mathbh{1}_{\{n_\gamma-1>S_{\gamma\{k,\ell\}}\}}\biggr]\\[-0.5pt]
&&{} +(1-z_k)(1-z_\ell)\mathrm{E}\biggl[
\frac{(n_\gamma-S_{\gamma\{k,\ell\}})(n_\gamma
-1-S_{\gamma\{k,\ell\}})}
{(N_\gamma-S_{\gamma\{k,\ell\}})(N_\gamma-1-S_{\gamma\{
k,\ell\}})}
\mathbh{1}_{\{n_\gamma>S_{\gamma\{k,\ell\}}\}}\biggr]\\[-0.5pt]
&\to& z_kz_\ell+(1-z_k)(1-z_\ell)\biggl(\frac{\rho-F(\tau)}{1-F(\tau
)}\biggr)^2\mathbh{1}_{\{\rho>F(\tau)\}}\\[-0.5pt]
& &{}+\{z_k(1-z_\ell)+(1-z_k)z_\ell\}\frac{(\rho-F(\tau
))\mathbh{1}_{\{\rho>F(\tau)\}}}{1-F(\tau)},\\[-0.5pt]
c_\gamma(y_k,y_{\ell})&=&d_\gamma(y_k,y_{\ell})-m'_\gamma
(y_k,y_\ell
)m'_\gamma(y_\ell,y_k)=\mathrm{o}_\gamma(1),\vspace*{-1pt}
\end{eqnarray*}
and A\ref{A3} is verified.

Cut-off sampling for $y_k\le\tau$ is essentially the complement of
stratified sampling with a ``take-all stratum'':
$I_{\gamma k}=1$ for the set $\{k\in U_\gamma\dvt z_k=1\}$.
Take-all strata are common in practice, particularly for the
highly-skewed populations in which cut-off sampling is attractive.
Arguments nearly identical to those above can be used to establish
A\ref
{A3} in the take-all case.
This take-all stratified design is analogous to the well-known class of
\textit{case--control studies} in epidemiology. We specifically consider
prospective case--control studies (e.g., Mantel~\cite{Mantel1973},
Langholz and Goldstein~\cite{LangholzGoldstein2001}, Arratia \textit
{et al.}~\cite{ArratiaGoldsteinLangholz2005}),
in which the finite population of all disease cases and controls is
stratified, disease cases ($z_k=1$) are selected with probability one,
and controls ($z_k=0$) are selected with probability less than one.

\subsection{With-replacement sampling with probability proportional to size}

Let $\{n_\gamma\}$ be a non-random sequence of positive integers with
$n_\gamma<N_\gamma$ and suppose that $f$ has strictly positive support:
$\int\mathbh{1}_{(-\infty,0]}f \,\mathrm{d}\lambda=0$. Consider
the case of
with-replacement sampling of $n_\gamma$ draws, with probability of
selecting element $k$ on the $h$th draw equal $p_{\gamma k}\in[0,1]$,
$\sum_{k\in U_\gamma}p_{\gamma k}=1$. While $p_{\gamma k}$ could be
constructed in many ways, a case of particular interest is $p_{\gamma
k}\propto Y_k$. This design is usually not feasible in practice, but
statistical agencies often attempt to draw samples with probability
proportional to a size measure (p.p.s.) that is highly correlated with
$Y$. Such a design will be highly efficient for estimation of the
$Y$-total (indeed, a fixed-size p.p.s. design with probabilities
proportional to $Y_k$ would exactly reproduce the $Y$-total).

For $h=1,\ldots,n_\gamma$, let $R_{\gamma h}$ be i.i.d. random
variables with
\[
\mathrm{P}(R_{\gamma h}=k\vert\mathcal{Y}_\gamma )=\frac{Y_k}{\sum
_{j\in U_\gamma}Y_j}.
\]
Then $I_{\gamma k}=\sum_{h=1}^{n_\gamma}\mathbh{1}_{\{R_{\gamma
h}=k\}
}$ counts the number of draws for which element $k$ is selected. Define
$W_{\gamma A}=N_\gamma^{-1}\sum_{j\in U_\gamma\dvtx j\notin A}Y_j$. Then
\begin{eqnarray*}
m_\gamma(y_k)&=&\frac{n_\gamma}{N_\gamma}y_k\mathrm{E}\biggl[\frac
{1}{N_\gamma
^{-1}y_k+W_{\gamma\{k\}}}\biggr],\\
m'_\gamma(y_k,y_\ell)&=&\frac{n_\gamma}{N_\gamma}y_k\mathrm
{E}\biggl[\frac{1}{N_\gamma
^{-1}(y_k+y_\ell)+W_{\gamma\{k,\ell\}}}\biggr],\\
v_\gamma(y_k)&=&\biggl(\frac{n_\gamma}{N_\gamma}y_k\biggr)^2\operatorname
{Var}\biggl(\frac
{1}{N_\gamma^{-1}y_k+W_{\gamma\{k\}}}\biggr)+\frac{n_\gamma}{N_\gamma}\frac{y_k}{N_\gamma}
\mathrm{E}\biggl[\frac{W_{\gamma\{k\}}}{(N_\gamma^{-1}y_k+W_{\gamma\{k\}}
)^2}\biggr],\\
c_\gamma(y_k,y_\ell)&=&\biggl(\frac{n_\gamma}{N_\gamma}\biggr)^2y_ky_\ell
\biggl\{
-\frac{1}{N_\gamma}\mathrm{E}\biggl[\frac{1}{(N_\gamma^{-1}(y_k+y_\ell
)+W_{\gamma
\{k,\ell\}})^2}\biggr]\\[-2pt]
&&\hspace*{57pt}{}+n_\gamma\operatorname{Var}\biggl(\frac{1}{ N_\gamma^{-1}(y_k+y_\ell
)+W_{\gamma\{k,\ell
\}}}\biggr)
\biggr\}.
\end{eqnarray*}

Under mild additional conditions, A\ref{A1} and A\ref{A2} can be
established using straightforward bounding and limiting arguments. A
sufficient set of conditions for either A\ref{A1} or A\ref{A2} is
$n_\gamma N_\gamma^{-1}\to\tau\in[0,1]$ as $\gamma\to\infty$ and
$\mathrm{E}[Y_1^6]<\infty$. Under these conditions, $m_\gamma
(y)=\tau y (\mathrm{E}[Y_1])^{-1}+\mathrm{o}_\gamma(1)$, and the
limiting c.d.f. is the same as in
length-biased sampling, as given by equation (\ref{eq:LBS}).\vspace*{-2pt}

\subsection{Endogenous stratification}\vspace*{-2pt}

Endogenous stratification, in which the sample is effectively
stratified on the value of the dependent variable, is common in the
health and social sciences (e.g., Hausman and Wise~\cite{HausmanWise1981}, Jewell~\cite{Jewell1985},
Shaw~\cite{Shaw1988}). Often, it arises by
design when a screening sample is
selected, the dependent variable is observed, and then covariates are
measured for a sub-sample that is stratified on the dependent
variable: for example, undersampling the high-income stratum
(Hausman and Wise~\cite{HausmanWise1981}). It can also arise through
uncontrolled selection
effects, in much the same way as length-biased sampling. One such
example comes from fisheries surveys, in which a field interviewer is
stationed at a dock for a fixed length of time, and intercepts
recreational fishing boats as they return to the dock. The interviewer
tends to select high-catch boats and, while busy measuring the fish
caught on those boats, misses more of the low-catch boats. Thus,
sampling effort is endogenously stratified on catch (Sullivan \textit
{et~al.}~\cite{Sullivan2006}).

We consider a sample endogenously stratified on the order statistics
of $Y$. Let $\{H_\gamma\}$ be a non-random sequence of positive
integers, which may remain bounded or go to infinity. For each~$\gamma
$, let $\{N_{\gamma h}\}_{h=1}^{H_\gamma}$ be a set of non-random
positive integers with $\sum_{h=1}^{H_\gamma}N_{\gamma h}=N_\gamma$,
and let $\{n_{\gamma h}\}_{h=1}^{H_\gamma}$ be a set of non-random
positive integers with $n_{\gamma h}\le N_{\gamma h}$. Let
\[
Y_{(1)}< Y_{(2)}< \cdots< Y_{(N_\gamma)}
\]
denote the order statistics for the $\gamma$th population, which is
stratified by taking the first $N_{\gamma1}$ values as stratum 1, the
next $N_{\gamma2}$ as stratum 2, etc., with the last $N_{\gamma
H_\gamma}$ values constituting stratum $H_\gamma$. The $\gamma$th
sample is then a stratified simple random sample without replacement of
size $n_{\gamma h}$ from the $N_{\gamma h}$ elements in stratum $h$.

Define $M_{\gamma0}=0$ and $M_{\gamma h}=\sum_{g=1}^h N_{\gamma h}$.
Because $H_\gamma,\ N_\gamma$ and $n_\gamma$ are not random, we then have
\begin{eqnarray*}
m_\gamma(y)&=&\sum_{h=1}^{H_\gamma} \frac{n_{\gamma h}}{N_{\gamma
h}}\mathrm{P}\bigl( Y_{(M_{\gamma,h-1})}<Y_k\le Y_{(M_{\gamma h})}\vert
Y_{k}=y\bigr)\\[-2pt]
&=& \sum_{h=1}^{H_\gamma} \frac{n_{\gamma h}}{N_{\gamma h}}
\mathrm{P}\biggl(\frac{M_{\gamma, h-1}}{N_{\gamma}-1}< F_{N_\gamma
-1}(y)\le\frac
{M_{\gamma, h}}{N_{\gamma}-1} \biggr),\vadjust{\goodbreak}
\end{eqnarray*}
where $F_{N_\gamma-1}(\cdot)$ is the empirical cumulative distribution
function for $ $ $\{Y_j\}_{j\in U_\gamma\dvtx j\ne k }$. From the
classical Glivenko--Cantelli theorem, $F_{N_\gamma-1}(y)$ converges
uniformly almost surely to~$F$. Similar computations can be used to
derive $m'_\gamma(y_1,y_2)$ and $c_\gamma(y_1,y_2)$ and their limits.
With such derivations, it is possible to establish the following
result, the proof of which is omitted.

\begin{result}\label{T5}
If $G(\alpha)=\lim_{\gamma\to\infty} \sum_{h=1}^{H_\gamma
}n_{\gamma
h}N_{\gamma h}^{-1}
\mathbh{1}_{(
N_{\gamma}^{-1}M_{\gamma,h-1},
N_{\gamma}^{-1}M_{\gamma h})}(\alpha)$
exists except for a finite number of points and is a piecewise
continuous nonnull function,
and the convergence is uniform in $\alpha$
then \textup{A\ref{A3}} and \textup{A\ref{A2}} hold.
\end{result}

\section{Conclusion}

We have given assumptions on the selection mechanism and the
superpopulation model under which the unweighted empirical c.d.f.
converges uniformly to a weighted version of the superpopulation c.d.f.
Because the conditions we specify on the informative selection
mechanism are closely tied to first and second-order inclusion
probabilities in a standard design-based survey sampling setting, the
conditions are verifiable. Our examples illustrate the computations for
selection mechanisms encountered in real surveys and observational
studies. We expect these conditions to be useful in studying the
properties of other basic sample statistics under informative
selection, which will be the subject of further research.

\begin{appendix}

\section{\texorpdfstring{Proofs of Theorems \protect\ref{T1} and \protect\ref{T2}}
{Proofs of Theorems 1 and 2}}\label{appa}

The first subsection contains the proof of Theorem~\ref{T1}.
The proof consists in showing the uniform $L_2$ convergence of the
empirical c.d.f., seen as a ratio of two random variables.
First, we show that from A\ref{A1} we can deduce the $L_2$ convergence
of both the numerator and denominator,
then the classical proof of Glivenko--Cantelli is adapted to
obtain a uniform $L_2$ convergence.

The second subsection contains the proof of Theorem~\ref{T2}.
We first construct two sequences of random variables $(\mathcal
{I}_\gamma ')$ and $Y'$
such that
$\forall\gamma, (\mathcal{I}_\gamma ', \mathcal{Y}_\gamma ')$ and
$(\mathcal{I}_\gamma , \mathcal{Y}_\gamma )$ have the same
distribution.
We then prove uniform $L_2$ convergence of the empirical c.d.f. defined
from $(\mathcal{I}_\gamma ')$ and $Y'$, almost surely in $Y'$. The
result is
``design-based'' in the sense that it is conditional on $Y'$, and is of
independent interest. We conclude by showing the almost sure convergence.

\subsection{\texorpdfstring{Proof of Theorem \protect\ref{T1}: Uniform $L_2$ convergence of the empirical c.d.f.}
{Proof of Theorem 1: Uniform L2 convergence of the empirical c.d.f.}}


\begin{lemma} \label{L1} Given a bounded measurable function
$b\dvtx\mathbb{R}\to\mathbb{R}$, \textup{A\ref{A0}} and \textup{A\ref{A1}} imply that
\[
\frac{\sum_{k\in U_\gamma}b(Y_k )I_{\gamma k}}{N_\gamma}
\mathop{\rightarrow}\limits_{\gamma\to\infty}^{{L_2}}
\int bmf \,\mathrm{d}\lambda.\vadjust{\goodbreak}
\]
\end{lemma}

\begin{pf}
Assume A\ref{A0} and A\ref{A1}.
The exchangeability property (\ref{eq:exchangeability}) implies
\[
\mathrm{E}\biggl[\frac{\sum_{k\in U_\gamma}b(Y_k )I_{\gamma k}}{N_\gamma}\biggr]=
\frac{\sum_{k\in U_\gamma} \mathrm{E}[b(Y_k )I_{\gamma
k}]}{N_\gamma}
= \int b m_\gamma f \,\mathrm{d}\lambda
\mathop{\rightarrow}\limits_{\gamma\to\infty}\int b m f \,\mathrm
{d}\lambda
\]
by (\ref{A0.1}), (\ref{A0.2}) and the dominated convergence theorem.
Further, (\ref{eq:exchangeability}) implies
\begin{eqnarray*}
&&\operatorname{Var}\biggl(\frac{\sum_{k\in U_\gamma}b(Y_k
)I_{\gamma k}}{N_\gamma
}\biggr)\\
&&\quad=\frac{1}{N_\gamma^2}\sum_{k,\ell\in U_\gamma}
\{
\operatorname{Cov}(b(Y_k )\mathrm{E}[I_{\gamma k}\vert Y_k , Y_{\ell}],
b(Y_{\ell})\mathrm{E}[I_{\gamma l}\vert Y_k , Y_{\ell}])
\\
& &\hspace*{44pt}\qquad{} + \mathrm{E}[b(Y_k )b(Y_{\ell})\operatorname{Cov}({I_{\gamma k},
I_{\gamma\ell}\vert Y_k , Y_{\ell}})] \}\\
&&\quad=\biggl(1-\frac{1}{N_\gamma}\biggr)
\biggl(\int b(y_1)b(y_2)m_\gamma'(y_1,y_2)m_\gamma
'(y_2,y_1)f(y_1)f(y_2)\,\mathrm{d}y_1 \,\mathrm{d}y_2\\
&& \hspace*{52pt}\qquad{}-\biggl(\int b(y_1)m_\gamma
'(y_1,y_2)f(y_1)f(y_2)\,\mathrm{d}y_1 \,\mathrm{d}y_2\biggr)^2\\
&&\hspace*{72pt}{} +\int b(y_1)b(y_2)c_\gamma
(y_1,y_2)f(y_1)f(y_2)\,\mathrm{d}y_1 \,\mathrm{d}y_2\biggr)\\
&&\qquad{}+\frac{1}{N_\gamma}
\biggl(\int b^2v_\gamma f \,\mathrm{d}\lambda+
\int b^2m_\gamma^2f \,\mathrm{d}\lambda-
\biggl(\int bm_\gamma f \,\mathrm{d}\lambda\biggr)^2\biggr)\\
&&\quad=\biggl(1-\frac{1}{N_\gamma}\biggr)
\biggl(\int b(y_1)b(y_2)\bigl(m_\gamma'(y_1,y_2)m_\gamma'(y_2,y_1)\\
&&\hspace*{110pt}\qquad {}-m_\gamma(y_1)m_\gamma
(y_2)\bigr)f(y_1)f(y_2)\,\mathrm{d}y_1 \,\mathrm{d}y_2\\
&&\hspace*{50pt}\qquad{} +\int b(y_1)b(y_2)c_\gamma
(y_1,y_2)f(y_1)f(y_2)\,\mathrm{d}y_1 \,\mathrm{d}y_2\biggr)\\
&&\qquad{}+\frac{1}{N_\gamma}
\biggl(\int b^2(v_\gamma+m_\gamma^2)f \,\mathrm{d}\lambda
-\biggl(\int bm_\gamma f \,\mathrm{d}\lambda\biggr)^2\biggr)\\
&&\quad=\mathrm{o}_\gamma(1)
\end{eqnarray*}
by (\ref{A1.1}), (\ref{A1.2}), and (\ref{A1.3}), and the result is proved.
\end{pf}


\begin{lemma}\label{L2} Under \textup{A\ref{A0}} and \textup{A\ref{A1}}, the numerator of
the empirical c.d.f. converges uniformly in $L_2$:
\[
\lim_{\gamma\to\infty}
\mathrm{E}\biggl[\biggl(\sup_{\alpha\in\mathbb{R}}
 \biggl|\frac{\sum_{k\in U_\gamma}\mathbh{1}_{(-\infty,\alpha]}(Y_k
)I_{\gamma k}}
{N_\gamma}-\int\mathbh{1}_{(-\infty,\alpha]} m_\gamma f \,\mathrm
{d}\lambda
\biggr|\biggr)^2\biggr]=0.
\]
\end{lemma}

\begin{pf}
We first define $G_\gamma\dvtx \mathbb{R}\to\mathbb{R}^+$ and
$G_s\dvtx \mathbb
{R}\to\mathbb{R}^+$
as
\[
G_\gamma(\alpha)=\frac{1}{N_\gamma}\sum_{k\in U_\gamma}\mathbh
{1}_{(-\infty,\alpha]}(Y_k )I_{\gamma k} \quad\mbox{and}\quad
G_s(\alpha)=\int\mathbh{1}_{(-\infty,\alpha]} mf \,\mathrm
{d}\lambda.
\]
With these definitions,
\[
\sup_{\alpha\in\mathbb{R}}
\biggl|\frac{\sum_{k\in U_\gamma}\mathbh{1}_{(-\infty,\alpha]}(Y_k
)I_{\gamma k}}
{N_\gamma}-\int\mathbh{1}_{(-\infty,\alpha]} m_\gamma f \,\mathrm
{d}\lambda
\biggr|=\|G_\gamma-G_s\|_\infty.
\]

Let $\eta\in\mathbb{N}^*$ index the positive integers and define a
sequence of subdivisions $\{\alpha_{\eta, q}\}_{q=0}^{\eta+1}$ of
$\mathbb{R}$ via $\alpha_{\eta, 0}=-\infty$, $\alpha_{\eta, \eta
+1}=\infty$, and for $q=1,\ldots,\eta$,
\[
\alpha_{\eta,q}=\inf\{\alpha\in\mathbb{R}|
G_s(\alpha)\geq\eta^{-1}{q G_s(\infty)}\}.
\]

We first show that for all positive integers $\eta$,
\[
\sup_{\alpha\in\mathbb{R}}\{|G_\gamma(\alpha)-G_s(\alpha)
|\}
\leq\max_{0\le q\le\eta+1}
\{|G_\gamma(\alpha_{\eta,q})-G_s(\alpha_{\eta,q})
|\}
+\frac{G_s(\infty)}{\eta}.
\]

Let $\eta\in\mathbb{N}$ and $\alpha\in\mathbb{R}$.
Then $\alpha\in[\alpha_{\eta,q},\alpha_{\eta,q+1}]$ for some
$0\le
q\le\eta$, and
\begin{eqnarray*}
G_\gamma(\alpha_{\eta,q}) &\leq& G_\gamma(\alpha) \leq
G_\gamma
(\alpha_{\eta,q+1}),\\
G_s(\alpha_{\eta,q}) &\leq& G_s(\alpha) \leq G_s(\alpha_{\eta
,q+1}),\\
G_s(\alpha_{\eta,q+1})-\frac{G_s(\infty)}{\eta} &\leq& G_s(\alpha)
\leq G_s(\alpha_{\eta,q})+\frac{G_s(\infty)}{\eta}.
\end{eqnarray*}
Combining these inequalities, we have
\begin{eqnarray*}
G_\gamma(\alpha_{\eta,q})-G_s(\alpha_{\eta,q})-\frac{G_s(\infty
)}{\eta}
&\leq& G_\gamma(\alpha)-G_s(\alpha) \\
&\leq& G_\gamma(\alpha_{\eta,q+1})-G_s(\alpha_{\eta,q+1})+\frac
{G_s(\infty)}{\eta},
\end{eqnarray*}
so that
\begin{eqnarray*}
&&|G_\gamma(\alpha)-G_s(\alpha)|\\
&&\quad\leq \max\{|G_\gamma(\alpha_{\eta,q})-G_s(\alpha_{\eta
,q})|,
|G_\gamma(\alpha_{\eta,q+1})-G_s(\alpha_{\eta,q+1})|
\}+\frac{G_s(\infty)}{\eta}\\
&&\quad\leq \max_{0\le q'\le\eta+1}
\{|G_\gamma(\alpha_{\eta,q'})-G_s(\alpha_{\eta,q'})
|\}+\frac{G_s(\infty)}{\eta}.
\end{eqnarray*}
Thus, for all $\alpha\in\mathbb{R}$,
\[
|G_\gamma(\alpha)-G_s(\alpha)|^2
\leq2\biggl( \max_{0\le q'\le\eta+1}
\{|G_\gamma(\alpha_{\eta,q'})-G_s(\alpha_{\eta,q'})
|^2\}
+\frac{G_s(\infty)^2}{\eta^2}\biggr),
\]
so that
%
%
\renewcommand{\theequation}{\arabic{equation}}
\setcounter{equation}{8}
\begin{equation}
\mathrm{E}[\|G_\gamma-G_s \|^2_\infty]
\leq2\mathrm{E}\Bigl[\max_{0\le q\le\eta+1}
\{|G_\gamma(\alpha_{\eta,q})-G_s(\alpha_{\eta,q})
|^2\}\Bigr]
+\frac{2G_s(\infty)^2}{\eta^2}.
\label{eq:L2sup}
\end{equation}

Let $\varepsilon>0$ be given. Choose $\eta\in\mathbb{N}$ so large that
$2G_s(\infty)^2\eta^{-2}<\varepsilon/2$, then use Lemma~\ref{L1} to
choose $\Gamma$ so that $\gamma\ge\Gamma$ implies
\[
2\mathrm{E}\Bigl[\max_{0\le q\le\eta+1}
\{|G_\gamma(\alpha_{\eta,q})-G_s(\alpha_{\eta,q})
|^2\}\Bigr]<\frac{\varepsilon}{2}.
\]
Hence, for all $\gamma\ge\Gamma$, the right-hand side of (\ref
{eq:L2sup}) is bounded by $\varepsilon$, which was arbitrary, so
$\lim_{\gamma\to\infty} \mathrm{E}[(\|G_\gamma-G_s \|_\infty
)^2 ]=0$.
\end{pf}


\begin{pf*}{Proof of Theorem \protect\ref{T1}}
By Definitions~\ref{D1} and~\ref{D3} and A\ref{A0}, for all $\alpha
\in
\mathbb{R}$,
\[
F_\gamma(\alpha) = \frac{G_\gamma(\alpha)}{G_\gamma(\infty
)+\mathbh
{1}_{G_\gamma(\infty)=0}}, \qquad
F_s(\alpha)=\frac{G_s(\alpha)}{G_s(\infty)},
\]
so
\begin{eqnarray*}
\|F_\gamma-F_s\|_\infty&=&\biggl\|\frac{G_\gamma
}{G_\gamma(\infty)+\mathbh{1}_{G_\gamma(\infty)=0}}-\frac
{G_s}{G_s(\infty)}\biggr\|_\infty\\
&= &\biggl\|\frac{G_\gamma-G_s}{G_s(\infty)}
+G_\gamma\frac{G_s(\infty)-(G_\gamma(\infty)+\mathbh{1}_{G_\gamma
(\infty)=0})}
{G_s(\infty)(G_\gamma(\infty)+\mathbh{1}_{G_\gamma(\infty)=0})}\biggr\|
_\infty\\
&\leq&\frac{\|G_\gamma-G_s\|_\infty}{G_s(\infty)}
+\frac{\|G_\gamma\|_\infty}{G_\gamma(\infty)+\mathbh
{1}_{G_\gamma(\infty)=0}}
\frac{|G_\gamma(\infty)+\mathbh{1}_{G_\gamma(\infty
)=0}-G_s(\infty
)|}{G_s(\infty)}\\
&\leq&\frac{\|G_\gamma-G_s\|_\infty}{G_s(\infty)}
+\frac{|G_\gamma(\infty)+\mathbh{1}_{G_\gamma(\infty
)=0}-G_s(\infty)|}{G_s(\infty)}\\
&\leq&\frac{\|G_\gamma-G_s\|_\infty}{G_s(\infty)}
+\frac{|G_s(\infty)-G_\gamma(\infty)|}{G_s(\infty)}
+\frac{\mathbh{1}_{G_\gamma(\infty)=0}}{G_s(\infty)}.
\end{eqnarray*}
From Lemma~\ref{L2}, the first two summands converge to 0 in $L_2$.
From (\ref{A1.5}), so does the third summand.
\end{pf*}

\subsection{\texorpdfstring{Proof of Theorem \protect\ref{T2}: Uniform almost sure convergence of the empirical c.d.f.}
{Proof of Theorem 2: Uniform almost sure convergence of the empirical c.d.f.}}

\subsubsection*{Construction of $\mathcal{I}_\gamma '$, $Y'$}

We define $Y'$ and $\mathcal{I}_\gamma '$ on the probability space
$(\Omega\times[0,1],\mathscr{A}\otimes\mathscr{B}_{[0,1]},\mathrm
{P}'=\mathrm{P}\otimes\lambda_{[0,1]})$.
First, define $Y'\dvtx\Omega\times[0,1]\to\mathbb{R}^{\mathbb
{N}}$ via
\[
Y'(\omega,x)=Y(\omega).
\]
Let $\mathcal{Y}_\gamma '$ be the vector of random variables
$(Y'_1,\ldots, Y'_{N_\gamma
})$ and note that
$\mathcal{Y}_\gamma '(\omega,x)=\mathcal{Y}_\gamma (\omega)$. Let
$S_{\gamma y}= \{i\in\mathbb
{N}^{N_\gamma}\dvt g_\gamma(i,y)\neq0\}$ and note that for a given
$y\in\mathbb{R}^{N_\gamma}$, $\sum_{i\in S_{\gamma y}}g_{\gamma}(i,y)=1$.
Define $h_\gamma\dvtx\mathbb{R}^{N_\gamma}\times\mathbb
{N}^{N_\gamma}\to
\mathbb{R}$ via
\[
h_\gamma(y,i)=\sup_{\alpha\in\mathbb{R}}
\biggl|\frac{\sum_{k\in U_\gamma}i_k \mathbh{1}_{(-\infty,\alpha]}(y_k)}
{\mathbh{1}_{i=0}+\sum_{k\in U_\gamma}(i_k)}-G_s(\alpha)\biggr|.
\]
We now impose an order on the $M_{\gamma y}$ vectors in $S_{\gamma y}$
by requiring $h_\gamma$
to be non-increasing; that is, for vectors $i^{(t)}, i^{(u)}\in
S_{\gamma y}$, $t<u$ if and only if
$h_\gamma(y,i^{(t)})\ge h_\gamma(y,i^{(u)})$. Any ties can be
resolved, for example, by randomization.
For $\omega\in\Omega$ and $x\in[0,1]$, we then define $\mathcal
{I}_\gamma '(\omega
,0)=i^{(1)}$ and for $x>0$
\[
\mathcal{I}_\gamma '(\omega,x)=
\sum_{u=1}^{M_{\gamma y}} i^{(u)}
\mathbh{1}_{(\sum_{t<u}g_\gamma(i^{(t)},\mathcal{Y}_\gamma (\omega
)),\sum_{t\leq
u}g_\gamma(i^{(t)},\mathcal{Y}_\gamma (\omega))]}(x).
\]
Because we use uniform measure on $\mathscr{B}_{[0,1]}$, the vector
$i^{(u)}$ is sampled from $S_{\gamma\mathcal{Y}_\gamma (\omega)}$
with probability
$g_\gamma(i^{(u)},\mathcal{Y}_\gamma (\omega))$.
Thus, by construction we have for all $\gamma$,
\[
\mathrm{P}'[\mathcal{I}_\gamma '=i\vert\mathcal{Y}_\gamma '
=y]=g_\gamma(i,y)= \mathrm{P}
[\mathcal{I}_\gamma =i\vert\mathcal{Y}_\gamma =y]
\]
and
$\mathrm{P}'[\mathcal{Y}_\gamma ' =y]=\mathrm{P}[\mathcal
{Y}_\gamma =y]$, so that
\[
\mathrm{P}'[\mathcal{I}_\gamma '=i, \mathcal{Y}_\gamma '
=y]=\mathrm{P}[\mathcal{I}_\gamma =i,\mathcal{Y}_\gamma =y].
\]
This yields the following property.
\begin{property}
For all $\gamma$,
\[
h_\gamma(\mathcal{Y}_\gamma ', \mathcal{I}_\gamma ')=\sup_{\alpha
\in\mathbb{R}}
|F'_\gamma(\alpha)-F_s(\alpha)|=\|F'_\gamma-F_s\|_\infty
\]
has the same law as $\|F_\gamma-F_s \|_\infty$, where $F'_\gamma$ is
defined in \textup{(\ref{eq:Fprime})}.
\end{property}

Define $G'_\gamma\dvtx\mathbb{R}\to\mathbb{R}^+$ via
\[
G'_\gamma(\alpha)=\frac{\sum_{k\in U_\gamma}\mathbh{1}_{(-\infty
,\alpha
]}(Y'_k)I'_{\gamma k}}{N_\gamma},
\]
noting that $F'_\gamma=G'_\gamma(G'_\gamma(\infty)+\mathbh
{1}_{G'_\gamma(\infty)=0})^{-1}$. We then have the following lemma.

\begin{lemma}\label{L3} Under \textup{A\ref{A0}} and \textup{A\ref{A2}}, for all
$\alpha
\in\mathbb{R}$,
\[
\lim_{\gamma\to\infty} \int_{[0,1]}
\bigl(G'_\gamma(\alpha)(\omega,x)
-G_s(\alpha)\bigr)^2 \,\mathrm{d}\lambda(x)=0\qquad \mbox{P-a.s. $(\omega)$}.
\]
\end{lemma}

\begin{pf}
Let
\[
\Omega_{\mathrm{GC}} =\biggl \{\omega\in\Omega\dvt
\lim_{\gamma\to\infty} \sup_{\alpha\in\mathbb{R}} \biggl|N_\gamma
^{-1}\sum_{k\in U_\gamma}\mathbh{1}_{(-\infty,\alpha]}(Y_k)(\omega)-
\int\mathbh{1}_{(-\infty,\alpha]} f \,\mathrm{d}\lambda\biggr|=0\biggr\}.
\]

From the Glivenko--Cantelli theorem, $\mathrm{P}(\Omega_{\mathrm{GC}})=1$.
We will show that for all $ \omega\in\Omega_{\mathrm{GC}}$,
\[
\int_{[0,1]}
\bigl(G'_\gamma(\alpha)(\omega,x)-G_s(\alpha)\bigr)^2
\,\mathrm{d}\lambda(x)
=\mathrm{o}_\gamma(1).
\]

Let $\omega\in\Omega_{\mathrm{GC}}$.
We then have
\begin{eqnarray*}
&&\sqrt{\int_{[0,1]}\bigl(G'_\gamma(\alpha)(\omega
,x)-G_s(\alpha
)\bigr)^2 \,\mathrm{d}\lambda(x)}\\
&&\quad\leq\sqrt{\int_{[0,1]} \biggl( G'_\gamma(\alpha)(\omega,x)
- \frac{\sum_{k\in U_\gamma}\mathbh{1}_{(-\infty,\alpha
]}(Y_k(\omega
) ) \int_{[0,1]}
I'_{\gamma k}(\omega,u)\,\mathrm{d}\lambda(u)}{N_\gamma} \biggr)^2 \,
\mathrm{d}\lambda
(x)}\\
&&\qquad{} + \biggl|\frac{\sum_{k\in U_\gamma}
\mathbh{1}_{(-\infty,\alpha]}(Y_k(\omega))\int_{[0,1]} I'_{\gamma
k}(\omega,u)\,\mathrm{d}\lambda(u)}{N_\gamma}
-\frac{\sum_{k\in U_\gamma}\mathbh
{1}_{(-\infty,\alpha]}(Y_k(\omega))m_\gamma(Y_k(\omega))}{N_\gamma
}\biggr|\\
&&\qquad{} +\biggl |\frac{\sum_{k\in U_\gamma}\mathbh{1}_{(-\infty
,\alpha]}(Y_k(\omega))m_\gamma(Y_k(\omega))}{N_\gamma}
-\int\mathbh{1}_{(-\infty,\alpha]} m_\gamma f \,\mathrm{d}\lambda
\biggr|\\
&&\qquad{} + \biggl|\int\mathbh{1}_{(-\infty,\alpha]} m_\gamma f
\,\mathrm{d}\lambda
-\int\mathbh{1}_{(-\infty,\alpha]} m f \,\mathrm{d}\lambda\biggr|.
\end{eqnarray*}
The first term is the square root of
\[
\operatorname{Var}\bigl(G_\gamma'(\alpha)\vert\mathcal{Y}_\gamma
'=(Y_1(\omega),\ldots,Y_{N_\gamma}(\omega
))\bigr)=N_\gamma^{-2}\mathrm{o}_\gamma(N_\gamma^2)=\mathrm{o}_\gamma(1)
\]
by (\ref{A2.1}).
The second term is
\[
\biggl|\sum_{k\in U_\gamma}
\frac{\mathbh{1}_{(-\infty,\alpha]} (Y_k(\omega)
)}{N_\gamma}
\bigl(\mathrm{E}[I'_{\gamma k}\vert\mathcal{Y}_\gamma ' = (Y_1 (\omega
),\ldots,
Y_{N_\gamma} (\omega))]
- m_\gamma(Y_k(\omega) ) \bigr)\biggr|
= \mathrm{o}_\gamma(1)
\]
by (\ref{A2.2}).
The third term is $\mathrm{o}_\gamma(1)$ because the convergence of the
empirical measure given by A\ref{A2} implies the convergence of the
integral for all bounded random variables. Finally, the fourth term is
$\mathrm{o}_\gamma(1)$ by A\ref{A0} and the dominated convergence theorem.
\end{pf}

The following lemma has its own interest, yielding design-based
uniform $L_2$ convergence of the empirical c.d.f.

\begin{lemma}\label{L4} Under \textup{A\ref{A0}} and \textup{A\ref{A2}},
\[
\int(h_\gamma(\mathcal{Y}_\gamma '(\omega,x),\mathcal{I}_\gamma
'(\omega,x)))^2 \,\mathrm{d}\lambda(x)=
\mathrm{o}_\gamma(1)\qquad \mbox{P-a.s.\ $(\omega)$}.
\]
\end{lemma}

\begin{pf}
Starting from Lemma~\ref{L3} and adapting the proof of Lemma~\ref{L2},
we have that:
$\mathrm{A\ref{A2}}\Rightarrow\int(\|G_\gamma(\mathcal{Y}_\gamma '(\omega
,x),\mathcal{I}_\gamma '(\omega
,x))-G_s\|_\infty)^2\,\mathrm{d}\lambda(x)=\mathrm{o}_\gamma(1)$
P-a.s. $(\omega)$. We then adapt the end of the proof of Theorem~\ref
{T1} and get the result.
\end{pf}

\begin{definition} For $\omega\in\Omega$, $\gamma\in\mathbb{N}$
and all
$\varepsilon>0$, $a_{\varepsilon,\gamma,\omega}\in[0,1]$ is
defined as
\[
a_{\varepsilon,\gamma,\omega}=\int_{[0,1]} \mathbh{1}_{\{h_\gamma
(\mathcal{Y}_\gamma
',\mathcal{I}_\gamma ')(\omega,x)\geq\varepsilon\}}\,\mathrm
{d}\lambda(x)
=\lambda_{[0,1]}\bigl(\{h_\gamma(\mathcal{Y}_\gamma ',\mathcal
{I}_\gamma ')(\omega,\cdot)\geq\varepsilon
\}\bigr).
\]
\end{definition}

\begin{property}
For all $\varepsilon>0$,
\[
\limsup_{\gamma\to\infty}
\mathbh{1}_{\{h_\gamma(\mathcal{Y}_\gamma ',\mathcal{I}_\gamma
')(\omega,x)>\varepsilon\}}=\mathbh
{1}_{\{0\}}\qquad
\mbox{P-a.s. $(\omega)$}.
\]
\end{property}
\begin{pf}
First note that $\forall x\in[0,1]$,
$\mathbh{1}_{\{h_\gamma(\mathcal{Y}_\gamma ',\mathcal{I}_\gamma
')(\omega,x)>\varepsilon\}}=\mathbh
{1}_{(0,a_{\varepsilon,\gamma,\omega}]}(x)$,
because by construction of $\mathcal{I}_\gamma ', \mathcal
{Y}_\gamma '$, $\{x\in[0,1]\dvt h_\gamma(\mathcal{Y}_\gamma
',\mathcal{I}_\gamma ')(\omega,x)>\varepsilon\}$
is a subinterval of $[0,1]$ containing $0$
of measure $a_{\varepsilon,\gamma,\omega}$. Further,
$\forall x\in[0,1]$,
%
%
\renewcommand{\theequation}{\arabic{equation}}
\setcounter{equation}{9}
\begin{equation}\limsup_{\gamma\to\infty} \mathbh{1}_{\{h_\gamma
(\mathcal{Y}_\gamma
',\mathcal{I}_\gamma ')(\omega,x)>\varepsilon\}}=\mathbh
{1}_{[0,\limsup_{\gamma\to\infty
} a_{\varepsilon,\gamma,\omega}]}(x).\label{eq:limsuph}
\end{equation}
By Lemma~\ref{L4}, the random variable
\[
h_\gamma(\mathcal{Y}_\gamma ',\mathcal{I}_\gamma ')(\omega,\cdot
)\dvtx\bigl([0,1],\mathscr{B}_{[0,1]},\lambda
_{[0,1]}\bigr)\to\mathbb{R}
\]
converges in $L_2(\lambda)$ to $0$, $\mbox{P-a.s. }(\omega)$, hence it
also converges in probability to $0$, and so $\lim_{\gamma\to\infty}
a_{\varepsilon,\gamma,\omega}=0$. The result then follows from equation
(\ref{eq:limsuph}).
\end{pf}


\begin{pf*}{Proof of Theorem \protect\ref{T2}}
We want to show that
\[
\mathrm{A\ref{A0}, A\ref{A2}}\Rightarrow\|F'_\gamma-F_s\|_\infty\,\mathop{\to}^{\mathrm{a.s.}}\,0 \qquad\mbox{as $\gamma\to\infty$},
\]
which is equivalent to showing that
\[
\textup{A\ref{A0}, A\ref{A2}}\Rightarrow\mathrm{P}'\Bigl(\Bigl\{\lim_{\gamma\to
\infty} h_\gamma(\mathcal{Y}_\gamma ',\mathcal{I}_\gamma ')=0\Bigr\}\Bigr)=1.
\]

Assume A\ref{A0} and A\ref{A2}. We calculate:
\begin{eqnarray*}
&&\mathrm{P}'\Bigl(\Bigl\{\lim_{\gamma\to\infty} h_\gamma(\mathcal
{Y}_\gamma ',\mathcal{I}_\gamma
')=0\Bigr\}\Bigr)\\
&&\quad=\mathrm{P}'\biggl(\bigcap_{\varepsilon>0} \bigcup_\Gamma\bigcap_{\gamma
>\Gamma} \{
h_\gamma(\mathcal{Y}_\gamma ',\mathcal{I}_\gamma ')<\varepsilon\}
\biggr)\\
&&\quad=\lim_{\varepsilon\to0} \mathrm{P}'\biggl(\bigcup_\Gamma\bigcap_{\gamma
>\Gamma} \{
h_\gamma(\mathcal{Y}_\gamma ',\mathcal{I}_\gamma ')<\varepsilon\}
\biggr)\\
&&\quad=\lim_{\varepsilon\to0} 1-\mathrm{P}'\biggl(\bigcap_\Gamma\bigcup_{\gamma
>\Gamma}
\{h_\gamma(\mathcal{Y}_\gamma ',\mathcal{I}_\gamma ')\geq
\varepsilon\}\biggr)\\
&&\quad=1-\lim_{\varepsilon\to0} \int\limsup_{\gamma\to\infty}
\mathbh
{1}_{\{h_\gamma(\mathcal{Y}_\gamma ',\mathcal{I}_\gamma ')(\omega
,x)\geq\varepsilon\}} \,\mathrm{d}\mathrm
{P}'(\omega,x).
\end{eqnarray*}

Let $\varepsilon>0$. Applying Fubini's theorem,
\begin{eqnarray*}
&&\int\limsup_{\gamma\to\infty} \mathbh
{1}_{\{
h_\gamma(\mathcal{Y}_\gamma ',\mathcal{I}_\gamma ')(\omega,x)\geq
\varepsilon\}} \,\mathrm{d}\mathrm{P}'(\omega,x)\\
&&\quad=\int\biggl(\int\limsup_\gamma\mathbh{1}_{\{h_\gamma(\mathcal
{Y}_\gamma ',\mathcal{I}_\gamma
')(\omega,x)\geq\varepsilon\}}
\,\mathrm{d}\lambda_{[0,1]}(x)\biggr) \,\mathrm{d}\mathrm{P}(\omega).
\end{eqnarray*}

Since we have
$\limsup_{\gamma\to\infty} \mathbh{1}_{\{h_\gamma(\mathcal
{Y}_\gamma ',\mathcal{I}_\gamma ')(\omega
,x)\geq\varepsilon\}}=\mathbh{1}_{\{0\}}(x)$ P-a.s. $(\omega)$, we also
have for all $\varepsilon>0$ that
\[
\int\limsup_{\gamma\to\infty} \mathbh{1}_{\{h_\gamma(\mathcal
{Y}_\gamma ',\mathcal{I}_\gamma ')(\omega
,x)\geq\varepsilon\}} \,\mathrm{d}\lambda_{[0,1]}(x)
=\int_{[0,1]}\mathbh{1}_{\{0\}}(x) \,\mathrm{d}\lambda_{[0,1]}(x)
=0
\]
P-a.s. $(\omega)$. Thus,
\[
\mathrm{P}'\Bigl(\Bigl\{\lim_{\gamma\to\infty} h_\gamma(\mathcal
{Y}_\gamma ',\mathcal{I}_\gamma ')=0
\Bigr\}\Bigr)=1.
\]
\upqed\end{pf*}

\section{\texorpdfstring{Proof of Corollaries \protect\ref{C1}, \protect\ref{C2}}
{Proof of Corollaries 1, 2}}\label{appb}
We state the following lemma which is a consequence of a theorem due
to P\'olya (e.g., Serfling~\cite{Serfling80}, page~18). The proof is omitted.
\begin{lemma}\label{lemmaf}
Let $\{u_\gamma(\cdot)\}_{\gamma\in\mathbb{N}}$ be a sequence of
increasing step functions, $u_\gamma\dvt\mathbb{R}\to[0,1]$, that
converges pointwise to a
continuous increasing function $u\dvtx\mathbb{R}\to[0,1]$ with $\lim
_{y\to
-\infty} u(y)=0$, $\lim_{y\to\infty} u(y)=1$ and
$0<u(y_1)=u(y_2)<1\Rightarrow y_1=y_2$.
Define $q_\gamma(p)=\inf\{y\in\mathbb{R}\dvt u_\gamma(y)\geq p\}$,
$q(p)=\inf\{y\in\mathbb{R}\dvt u(y)\geq p\}$.
Then for all $K$ a compact subset of $(0,1)$,
$\lim_{\gamma\to\infty} \sup_{p\in K}\{q_\gamma(p)-q(p)\}=0.$
\end{lemma}
%
\subsection{\texorpdfstring{Proof of Corollary \protect\ref{C1}}{Proof of Corollary 1}}
\begin{pf}
As $m_\gamma f$ and $mf$ may have different supports, we extend the
definition of
$\xi_s$ by
\[
\forall p\in\mathbb{R}, \qquad \xi_s(p)=\inf\{y\in\mathbb{R}\dvt
F_s(y)\geq p\}.
\]
Let $K$ be a compact subset of $(0,1)$. Then
\[
\sup_{p\in K}| \xi_\gamma(p)-\xi_s(p)|\mathop{\rightarrow}\limits_{\gamma\to
\infty}^{{P}}0
\]
if from all subsequences one can extract a subsequence that converges
a.s. to $0$.
Let $\tau\dvt\mathbb{N}\to\mathbb{N}$ be a strictly increasing function.
If $\|F_\gamma-F_s\|_\infty{\to}^{L_2}0$
then $\|F_{\tau(\gamma)}-F_s\|_\infty{\to}^{L_2}0$
and $\|F_{\tau(\gamma)}-F_s\|_\infty{\to}^{P}0$. Then there
exists $\rho\dvt\mathbb{N}\to\mathbb{N}$ strictly increasing
such that $\|F_{\tau(\rho(\gamma))}-F_s\|_\infty{\to}^{\mathrm {a.s.}}0$
and by Lemma~\ref{lemmaf},
$\mathrm{P}(\lim_{\gamma\to\infty}\sup_{p\in K}|\xi_{\tau(\rho
(\gamma
))}(p)-\xi_s(p)|=0)=1.$

For the uniform $L_2$ convergence, let $p\in(0,1)$ and $\alpha\in
\mathbb{R}$. Then
$|F_\gamma(\alpha)-F_s(\alpha)|\le\|F_\gamma-F_s\|_\infty,$
so that
\begin{eqnarray*}
\{\alpha\in\mathbb{R}\dvt F_s(\alpha)\ge p+\|F_\gamma-F_s\|
_\infty\}
&\subset&
\{\alpha\in\mathbb{R}\dvt F_\gamma(\alpha)\ge p\}\\
&\subset& \{\alpha\in\mathbb{R}\dvt F_s(\alpha)\ge p-\|F_\gamma
-F_s\|
_\infty\},
\end{eqnarray*}
and
\begin{eqnarray*}
\inf\{\alpha\in\mathbb{R}\dvt F_s(\alpha)\ge p+\|F_\gamma-F_s\|
_\infty\}
&\ge&
\inf\{\alpha\in\mathbb{R}\dvtx F_\gamma(\alpha)\ge p\}\\
&\ge&\inf\{\alpha\in\mathbb{R}\dvt F_s(\alpha)\ge p-\|F_\gamma
-F_s\|
_\infty\}.
\end{eqnarray*}

Hence, $\forall p\in(0,1),$
$\xi_s(p+\|F_\gamma-F_s\|_\infty)\ge\xi_\gamma(p)\ge\xi
_s(p-\|F_\gamma-F_s\|_\infty).$

Further, $f$ has compact support by hypothesis, so there exists
$b>0$ such that the supports of $(m_\gamma f)_{\gamma\in\mathbb{N}}$
and $mf$ are included in
$[-b,b]$. So $\forall p\in(0,1), $ $\gamma\in\mathbb{N},$ $-b\leq
\xi
_\gamma(p)\leq b$,
$-b\leq\xi_s(p)\leq b$. By combining these three inequalities, we
have, $\forall p\in(0,1)$:
%
%
\renewcommand{\theequation}{\arabic{equation}}
\setcounter{equation}{10}
\begin{equation}
|\xi_s (p ) - \xi_\gamma(p ) | \le\min\{
b,\xi_s (p + \|F_\gamma- F_s \|_\infty) \}
- \max\{ -b,\xi_s (p - \|F_\gamma- F_s \|_\infty
) \} .
\label{eq:bound1}
\end{equation}

Since $K\subset(0,1)$ is compact, there exists $a\in(0,1)$ such that
$K\subset[a,1-a]$. With the assumed continuity of $F_s$,
we have that $\xi_s$ is uniformly continuous on any subinterval of
$[0,1]$ that does not contain zero.
Thus, for $\varepsilon>0$, there exists $\eta\in(0,a/2)$
such that $p\in K$ implies $|\xi_s(p+\eta)-\xi_s(p-\eta)|\le
\varepsilon
.$ If $\|F_\gamma-F_s\|_\infty\le\eta$, then
$p+\|F_\gamma-F_s\|_\infty\le p+\eta<1-a/2$, and $\xi_s(p+\|
F_\gamma
-F_s\|_\infty)<b$,
$p-\|F_\gamma-F_s\|_\infty\ge p-\eta>a/2$ and $\xi_s(p-\|F_\gamma
-F_s\|
_\infty)>-b$,
so equation (\ref{eq:bound1}) is bounded by $\varepsilon$.
If $\|F_\gamma-F_s\|_\infty> \eta$, then (\ref{eq:bound1}) is bounded
by $(2b)$ $\mathbh{1}_{\{\|F_\gamma-F_s\|_\infty> \eta\}}$. Thus,
\[
\mathrm{E}\Bigl[\Bigl(\sup_{p\in K}| \xi_\gamma(p)-\xi_s(p)|\Bigr)^2\Bigr]\le
\varepsilon^2+4b^2\mathrm{P}(\|F_\gamma-F_s\|_\infty> \eta).
\]
Since $\varepsilon$ was arbitrary and $\mathrm{P}(\|F_\gamma-F_s\|
_\infty> \eta
)\to0$ as $\gamma\to\infty$, the result follows.
\end{pf}

\subsection{\texorpdfstring{Proof of Corollary \protect\ref{C2}}{Proof of Corollary 2}}

\begin{pf}
If $\|F'_\gamma-F_s\|_\infty{\to}^{\mathrm{a.s.}}0$, then
for all
$K$ a compact subset of $(0,1)$, and all $(\omega,x) \in
\{(\omega,x)\dvt\|F'_\gamma-F_s\|_\infty\to0\}$, we apply Lemma
\ref{lemmaf} with $u_\gamma=F_\gamma'(\omega,x)$, $u=F_s$, and
obtain that
$\mathrm{P}'(\lim_{\gamma\to\infty}\sup_{p\in K}|\xi'_\gamma
(p)-\xi'_s(p)|=0)=1$.
\end{pf}
\end{appendix}
\section*{Acknowledgement}
This research was supported in part by the US National Science Foundation (SES--0922142).

%

\printhistory


\begin{thebibliography}{40}

\bibitem{ArratiaGoldsteinLangholz2005}
%
\begin{barticle}[mr]
\bauthor{\bsnm{Arratia},~\bfnm{Richard}\binits{R.}},
\bauthor{\bsnm{Goldstein},~\bfnm{Larry}\binits{L.}} \AND
\bauthor{\bsnm{Langholz},~\bfnm{Bryan}\binits{B.}}
(\byear{2005}).
\btitle{Local central limit theorems, the high-order correlations of rejective
sampling and logistic likelihood asymptotics}.
\bjournal{Ann. Statist.}
\bvolume{33}
\bpages{871--914}.
\bid{doi={10.1214/009053604000000706}, issn={0090-5364}, mr={2163162}}
\bptok{imsref}%
\end{barticle}
%
\endbibitem

\bibitem{Binder1983}
%
\begin{barticle}[mr]
\bauthor{\bsnm{Binder},~\bfnm{David~A.}\binits{D.A.}}
(\byear{1983}).
\btitle{On the variances of asymptotically normal estimators from complex
surveys}.
\bjournal{Internat. Statist. Rev.}
\bvolume{51}
\bpages{279--292}.
\bid{doi={10.2307/1402588}, issn={0306-7734}, mr={0731144}}
\bptok{imsref}%
\end{barticle}
%
\endbibitem

\bibitem{Brecklingetal1994}
%
\begin{barticle}[auto:STB|2011/10/17|13:52:43]
\bauthor{\bsnm{Breckling},~\bfnm{J.~U.}\binits{J.U.}},
\bauthor{\bsnm{Chambers},~\bfnm{R.~L.}\binits{R.L.}},
\bauthor{\bsnm{Dorfman},~\bfnm{A.~H.}\binits{A.H.}},
\bauthor{\bsnm{Tam},~\bfnm{S.~M.}\binits{S.M.}} \AND
\bauthor{\bsnm{Welsh},~\bfnm{A.~H.}\binits{A.H.}}
(\byear{1994}).
\btitle{Maximum likelihood inference from sample survey data}.
\bjournal{Internat. Statist. Rev.}
\bvolume{62}
\bpages{349--363}.
\bptok{imsref}%
\end{barticle}
%
\endbibitem

\bibitem{bre00}
%
\begin{barticle}[mr]
\bauthor{\bsnm{Breidt},~\bfnm{F.~Jay}\binits{F.J.}} \AND
\bauthor{\bsnm{Opsomer},~\bfnm{Jean~D.}\binits{J.D.}}
(\byear{2000}).
\btitle{Local polynomial regresssion estimators in survey sampling}.
\bjournal{Ann. Statist.}
\bvolume{28}
\bpages{1026--1053}.
\bid{doi={10.1214/aos/1015956706}, issn={0090-5364}, mr={1810918}}
\bptok{imsref}%
\end{barticle}
%
\endbibitem

\bibitem{bre08b}
%
\begin{barticle}[mr]
\bauthor{\bsnm{Breidt},~\bfnm{F.~Jay}\binits{F.J.}} \AND
\bauthor{\bsnm{Opsomer},~\bfnm{Jean~D.}\binits{J.D.}}
(\byear{2008}).
\btitle{Endogenous post-stratification in surveys: Classifying with a
sample-fitted model}.
\bjournal{Ann. Statist.}
\bvolume{36}
\bpages{403--427}.
\bid{doi={10.1214/009053607000000703}, issn={0090-5364}, mr={2387977}}
\bptok{imsref}%
\end{barticle}
%
\endbibitem

\bibitem{CasselSarndalWretman1977}
%
\begin{bbook}[mr]
\bauthor{\bsnm{Cassel},~\bfnm{Claes-Magnus}\binits{C.M.}},
\bauthor{\bsnm{S{\"a}rndal},~\bfnm{Carl-Erik}\binits{C.E.}} \AND
\bauthor{\bsnm{Wretman},~\bfnm{Jan~H{\. a}kan}\binits{J.H.}}
(\byear{1977}).
\btitle{Foundations of Inference in Survey Sampling}.
\bseries{Wiley Series in Probability and Mathematical Statistics}.
\baddress{New York}: \bpublisher{Wiley-Interscience}.
\bid{mr={0652527}}
\bptok{imsref}%
\end{bbook}
%
\endbibitem

\bibitem{MR1616053}
%
\begin{barticle}[mr]
\bauthor{\bsnm{Chambers},~\bfnm{Raymond~L.}\binits{R.L.}},
\bauthor{\bsnm{Dorfman},~\bfnm{Alan~H.}\binits{A.H.}} \AND
\bauthor{\bsnm{Wang},~\bfnm{Suojin}\binits{S.}}
(\byear{1998}).
\btitle{Limited information likelihood analysis of survey data}.
\bjournal{J. R. Stat. Soc. Ser. B Stat. Methodol.}
\bvolume{60}
\bpages{397--411}.
\bid{doi={10.1111/1467-9868.00132}, issn={1369-7412}, mr={1616053}}
\bptok{imsref}%
\end{barticle}
%
\endbibitem

\bibitem{ChambersSkinner2003}
%
\begin{bbook}[mr]
\bauthor{\bsnm{Chambers},~\bfnm{Raymond~L.}\binits{R.L.}} \AND
\bauthor{\bsnm{Skinner},~\bfnm{C. J.}\binits{C. J.}}
(\byear{2003}).
\btitle{Analysis of Survey Data}.
\baddress{Chichester}: \bpublisher{Wiley}.
\bid{doi={10.1002/0470867205}, mr={1978840}}
\bptok{imsref}%
\end{bbook}
%
\endbibitem

\bibitem{cox69}
%
\begin{bincollection}[auto:STB|2011/10/17|13:52:43]
\bauthor{\bsnm{Cox},~\bfnm{D.}\binits{D.}}
(\byear{1969}).
\btitle{Some sampling problems in technology}.
In \bbooktitle{New Developments in Survey Sampling}
(\beditor{\bfnm{U.}\binits{U.}~\bsnm{Johnson}} \AND
\beditor{\bfnm{H.}\binits{H.}~\bsnm{Smith}}, eds.)
\bpages{506--527}.
\baddress{New York}: \bpublisher{Wiley Interscience}.
\bptok{imsref}%
\end{bincollection}
%
\endbibitem

\bibitem{EaglesonWeber1978}
%
\begin{barticle}[mr]
\bauthor{\bsnm{Eagleson},~\bfnm{G.~K.}\binits{G.K.}} \AND
\bauthor{\bsnm{Weber},~\bfnm{N.~C.}\binits{N.C.}}
(\byear{1978}).
\btitle{Limit theorems for weakly exchangeable arrays}.
\bjournal{Math. Proc. Cambridge Philos. Soc.}
\bvolume{84}
\bpages{123--130}.
\bid{issn={0305-0041}, mr={0501275}}
\bptok{imsref}%
\end{barticle}
%
\endbibitem

\bibitem{EidehNathan2009}
%
\begin{barticle}[mr]
\bauthor{\bsnm{Eideh},~\bfnm{Abdulhakeem}\binits{A.}} \AND
\bauthor{\bsnm{Nathan},~\bfnm{Gad}\binits{G.}}
(\byear{2009}).
\btitle{Two-stage informative cluster sampling-estimation and
prediction with
applications for small-area models}.
\bjournal{J. Statist. Plann. Inference}
\bvolume{139}
\bpages{3088--3101}.
\bid{doi={10.1016/j.jspi.2009.02.019}, issn={0378-3758}, mr={2535185}}
\bptok{imsref}%
\end{barticle}
%
\endbibitem

\bibitem{EidehNathan2007}
%
\begin{barticle}[mr]
\bauthor{\bsnm{Eideh},~\bfnm{Abdulhakeem A.~H.}\binits{A.A.H.}}
\AND
\bauthor{\bsnm{Nathan},~\bfnm{Gad}\binits{G.}}
(\byear{2006}).
\btitle{Fitting time series models for longitudinal survey data under
informative sampling}.
\bjournal{J. Statist. Plann. Inference}
\bvolume{136}
\bpages{3052--3069}.
\bid{doi={10.1016/j.jspi.2004.12.003}, issn={0378-3758}, mr={2256216}}
\bptok{imsref}%
\end{barticle}
%
\endbibitem

\bibitem{EidehNathan2007c}
%
\begin{bmisc}[mr]
\bauthor{\bsnm{Eideh},~\bfnm{Abdulhakeem A.~H.}\binits{A.A.H.}}
\AND
\bauthor{\bsnm{Nathan},~\bfnm{Gad}\binits{G.}}
(\byear{2007}).
\bhowpublished{Corrigendum to ``Fitting time series models for longitudinal survey data under
informative sampling'' [\textit{J. Statist. Plann. Inference} \textbf{136} 3052--3069].
\textit{J. Statist. Plann. Inference} \textbf{137} 628}.
\bptok{imsref}%
\end{bmisc}
%
\endbibitem

\bibitem{Fuller2009}
%
\begin{bbook}[auto:STB|2011/10/17|13:52:43]
\bauthor{\bsnm{Fuller},~\bfnm{W.}\binits{W.}}
(\byear{2009}).
\btitle{Sampling Statistics}.
\baddress{New York}: \bpublisher{Wiley}.
\bptok{imsref}%
\end{bbook}
%
\endbibitem

\bibitem{Hajek81}
%
\begin{bbook}[mr]
\bauthor{\bsnm{H{\'a}jek},~\bfnm{Jaroslav}\binits{J.}}
(\byear{1981}).
\btitle{Sampling from a Finite Population}.
\bseries{Statistics: Textbooks and Monographs}
\bvolume{37}.
\baddress{New York}: \bpublisher{Dekker}.
\bnote{Edited by V{\'a}clav Dupa{\v{c}}, With a foreword by P. K. Sen}.
\bid{mr={0627744}}
\bptok{imsref}%
\end{bbook}
%
\endbibitem

\bibitem{HausmanWise1981}
%
\begin{bincollection}[auto:STB|2011/10/17|13:52:43]
\bauthor{\bsnm{Hausman},~\bfnm{J.}\binits{J.}} \AND
\bauthor{\bsnm{Wise},~\bfnm{D.}\binits{D.}}
(\byear{1981}).
\btitle{Stratification on endogenous variables and estimation: The
Gary income
maintenance experiment}.
In \bbooktitle{Structural Analysis of Discrete Data}.
\baddress{Cambridge, MA}: \bpublisher{MIT Press}.
\bptok{imsref}%
\end{bincollection}
%
\endbibitem

\bibitem{IsakiFuller82}
%
\begin{barticle}[mr]
\bauthor{\bsnm{Isaki},~\bfnm{Cary~T.}\binits{C.T.}} \AND
\bauthor{\bsnm{Fuller},~\bfnm{Wayne~A.}\binits{W.A.}}
(\byear{1982}).
\btitle{Survey design under the regression superpopulation model}.
\bjournal{J.~Amer. Statist. Assoc.}
\bvolume{77}
\bpages{89--96}.
\bid{issn={0162-1459}, mr={0648029}}
\bptok{imsref}%
\end{barticle}
%
\endbibitem

\bibitem{Jewell1985}
%
\begin{barticle}[mr]
\bauthor{\bsnm{Jewell},~\bfnm{Nicholas~P.}\binits{N.P.}}
(\byear{1985}).
\btitle{Least squares regression with data arising from stratified
samples of
the dependent variable}.
\bjournal{Biometrika}
\bvolume{72}
\bpages{11--21}.
\bid{doi={10.1093/biomet/72.1.11}, issn={0006-3444}, mr={0790196}}
\bptok{imsref}%
\end{barticle}
%
\endbibitem

\bibitem{KishFrankel1974}
%
\begin{barticle}[mr]
\bauthor{\bsnm{Kish},~\bfnm{Leslie}\binits{L.}} \AND
\bauthor{\bsnm{Frankel},~\bfnm{Martin~Richard}\binits{M.R.}}
(\byear{1974}).
\btitle{Inference from complex samples}.
\bjournal{J. Roy. Statist. Soc. Ser. B}
\bvolume{36}
\bpages{1--37}.
\bid{issn={0035-9246}, mr={0365812}}
\bptnote{check related}
\bptok{imsref}%
\end{barticle}
%
\endbibitem

\bibitem{KriegerPfeffermann1992}
%
\begin{barticle}[auto:STB|2011/10/17|13:52:43]
\bauthor{\bsnm{Krieger},~\bfnm{A.~M.}\binits{A.M.}} \AND
\bauthor{\bsnm{Pfeffermann},~\bfnm{D.}\binits{D.}}
(\byear{1992}).
\btitle{Maximum likelihood estimation for complex sample surveys}.
\bjournal{Survey Methodology}
\bvolume{18}
\bpages{225--239}.
\bptok{imsref}%
\end{barticle}
%
\endbibitem

\bibitem{LangholzGoldstein2001}
%
\begin{barticle}[pbm]
\bauthor{\bsnm{Langholz},~\bfnm{B.}\binits{B.}} \AND
\bauthor{\bsnm{Goldstein},~\bfnm{L.}\binits{L.}}
(\byear{2001}).
\btitle{Conditional logistic analysis of case--control studies with complex
sampling}.
\bjournal{Biostatistics}
\bvolume{2}
\bpages{63--84}.
\bid{doi={10.1093/biostatistics/2.1.63}, issn={1468-4357}, pii={2/1/63}, pmid={12933557}}
\bptok{imsref}%
\end{barticle}
%
\endbibitem

\bibitem{Leigh1988}
%
\begin{barticle}[auto:STB|2011/10/17|13:52:43]
\bauthor{\bsnm{Leigh},~\bfnm{G.~M.}\binits{G.M.}}
(\byear{1988}).
\btitle{A comparison of estimates of natural mortality from fish tagging
experiments}.
\bjournal{Biometrika}
\bvolume{75}
\bpages{347--353}.
\bptok{imsref}%
\end{barticle}
%
\endbibitem

\bibitem{Mantel1973}
%
\begin{barticle}[pbm]
\bauthor{\bsnm{Mantel},~\bfnm{N.}\binits{N.}}
(\byear{1973}).
\btitle{Synthetic retrospective studies and related topics}.
\bjournal{Biometrics}
\bvolume{29}
\bpages{479--486}.
\bid{issn={0006-341X}, pmid={4793136}}
\bptok{imsref}%
\end{barticle}
%
\endbibitem

\bibitem{NowellStanley1991}
%
\begin{barticle}[auto:STB|2011/10/17|13:52:43]
\bauthor{\bsnm{Nowell},~\bfnm{C.}\binits{C.}} \AND
\bauthor{\bsnm{Stanley},~\bfnm{L.~R.}\binits{L.R.}}
(\byear{1991}).
\btitle{Length-biased sampling in mall intercept surveys}.
\bjournal{Journal of Marketing Research}
\bvolume{28}
\bpages{475--479}.
\bptok{imsref}%
\end{barticle}
%
\endbibitem

\bibitem{PatilRao1978}
%
\begin{barticle}[mr]
\bauthor{\bsnm{Patil},~\bfnm{G.~P.}\binits{G.P.}} \AND
\bauthor{\bsnm{Rao},~\bfnm{C.~R.}\binits{C.R.}}
(\byear{1978}).
\btitle{Weighted distributions and size-biased sampling with
applications to
wildlife populations and human families}.
\bjournal{Biometrics}
\bvolume{34}
\bpages{179--189}.
\bid{doi={10.2307/2530008}, issn={0006-341X}, mr={0507202}}
\bptok{imsref}%
\end{barticle}
%
\endbibitem

\bibitem{PfeffermannKriegerRinott98}
%
\begin{barticle}[mr]
\bauthor{\bsnm{Pfeffermann},~\bfnm{Danny}\binits{D.}},
\bauthor{\bsnm{Krieger},~\bfnm{Abba~M.}\binits{A.M.}} \AND
\bauthor{\bsnm{Rinott},~\bfnm{Yosef}\binits{Y.}}
(\byear{1998}).
\btitle{Parametric distributions of complex survey data under informative
probability sampling}.
\bjournal{Statist. Sinica}
\bvolume{8}
\bpages{1087--1114}.
\bid{issn={1017-0405}, mr={1666233}}
\bptok{imsref}%
\end{barticle}
%
\endbibitem

\bibitem{PfeffermanMouraDaSilvaLuisdonascimento2006}
%
\begin{barticle}[mr]
\bauthor{\bsnm{Pfeffermann},~\bfnm{Danny}\binits{D.}},
\bauthor{\bsnm{Moura},~\bfnm{Fernando Antonio Da~Silva}\binits{F.A.D.S.}}
\AND\bauthor{\bsnm{Silva},~\bfnm{Pedro Luis do~Nascimento}\binits
{P.L.d.N.}}
(\byear{2006}).
\btitle{Multi-level modelling under informative sampling}.
\bjournal{Biometrika}
\bvolume{93}
\bpages{943--959}.
\bid{doi={10.1093/biomet/93.4.943}, issn={0006-3444}, mr={2285081}}
\bptok{imsref}%
\end{barticle}
%
\endbibitem

\bibitem{PfeffermannSverchkov1999}
%
\begin{barticle}[mr]
\bauthor{\bsnm{Pfeffermann},~\bfnm{Danny}\binits{D.}} \AND
\bauthor{\bsnm{Sverchkov},~\bfnm{Michail}\binits{M.}}
(\byear{1999}).
\btitle{Parametric and semi-parametric estimation of regression models fitted
to survey data}.
\bjournal{Sankhy\=a Ser. B}
\bvolume{61}
\bpages{166--186}.
\bid{issn={0581-5738}, mr={1720710}}
\bptok{imsref}%
\end{barticle}
%
\endbibitem

\bibitem{PfeffermannSverchkov2007}
%
\begin{barticle}[mr]
\bauthor{\bsnm{Pfeffermann},~\bfnm{Danny}\binits{D.}} \AND
\bauthor{\bsnm{Sverchkov},~\bfnm{Michail}\binits{M.}}
(\byear{2007}).
\btitle{Small-area estimation under informative probability sampling
of areas
and within the selected areas}.
\bjournal{J. Amer. Statist. Assoc.}
\bvolume{102}
\bpages{1427--1439}.
\bid{doi={10.1198/016214507000001094}, issn={0162-1459}, mr={2412558}}
\bptok{imsref}%
\end{barticle}
%
\endbibitem

\bibitem{PfeffermannSverchkov2009}
%
\begin{bincollection}[auto:STB|2011/10/17|13:52:43]
\bauthor{\bsnm{Pfeffermann},~\bfnm{D.}\binits{D.}} \AND
\bauthor{\bsnm{Sverchkov},~\bfnm{M.}\binits{M.}}
(\byear{2009}).
\btitle{Inference under informative sampling}.
In \bbooktitle{Handbook of Statistics}
(\beditor{\bfnm{D.}\binits{D.}~\bsnm{Pfefferman}} \AND
\beditor{\bfnm{C.~R.}\binits{C.R.}~\bsnm{Rao}}, eds.)
\bvolume{29B}
\bpages{455--487}.
\baddress{Amsterdam}: \bpublisher{North-Holland}.
\bptok{imsref}%
\end{bincollection}
%
\endbibitem

\bibitem{PfeffermannSverchkov2003}
%
\begin{bincollection}[mr]
\bauthor{\bsnm{Pfeffermann},~\bfnm{Danny}\binits{D.}} \AND
\bauthor{\bsnm{Sverchkov},~\bfnm{M.~Yu.}\binits{M.Y.}}
(\byear{2003}).
\btitle{Fitting generalized linear models under informative sampling}.
In \bbooktitle{Analysis of Survey Data}.
\bseries{Wiley Ser. Surv. Methodol.}
\bpages{175--195}.
\baddress{Chichester}: \bpublisher{Wiley}.
\bid{doi={10.1002/0470867205.ch12}, mr={1978851}}
\bptok{imsref}%
\end{bincollection}
%
\endbibitem

\bibitem{Poisson1837}
%
\begin{bincollection}[auto:STB|2011/10/17|13:52:43]
\bauthor{\bsnm{Poisson},~\bfnm{S.~D.}\binits{S.D.}}
(\byear{1837}).
\btitle{Recherches sur la probabilit\'e des jugements en mati\`ere criminelle
et en mati\`ere civile}.
In \bbooktitle{Proc\'ed\'es des R\`egles G\'en\'erales du Calcul des
Probabiliti\'es}.
\baddress{Paris}: \bpublisher{Bachelier, Imprimeur-Libraire pour les Math\'ematiques}.
\bptok{imsref}%
\end{bincollection}
%
\endbibitem

\bibitem{rob83}
%
\begin{barticle}[mr]
\bauthor{\bsnm{Robinson},~\bfnm{P.~M.}\binits{P.M.}} \AND
\bauthor{\bsnm{S{\"a}rndal},~\bfnm{Carl-Erik}\binits{C.E.}}
(\byear{1983}).
\btitle{Asymptotic properties of the generalized regression estimator in
probability sampling}.
\bjournal{Sankhy\=a Ser. B}
\bvolume{45}
\bpages{240--248}.
\bid{issn={0581-5738}, mr={0748468}}
\bptok{imsref}%
\end{barticle}
%
\endbibitem

\bibitem{SarndalSwenssonWretman92}
%
\begin{bbook}[mr]
\bauthor{\bsnm{S{\"a}rndal},~\bfnm{Carl-Erik}\binits{C.E.}},
\bauthor{\bsnm{Swensson},~\bfnm{Bengt}\binits{B.}} \AND
\bauthor{\bsnm{Wretman},~\bfnm{Jan}\binits{J.}}
(\byear{1992}).
\btitle{Model Assisted Survey Sampling}.
\baddress{New York}: \bpublisher{Springer}.
\bid{mr={1140409}}
\bptok{imsref}%
\end{bbook}
%
\endbibitem

\bibitem{Serfling80}
%
\begin{bbook}[mr]
\bauthor{\bsnm{Serfling},~\bfnm{Robert~J.}\binits{R.J.}}
(\byear{1980}).
\btitle{Approximation Theorems of Mathematical Statistics}.
\bseries{Wiley Series in Probability and Mathematical Statistics}.
\baddress{New York}: \bpublisher{Wiley}.
\bid{mr={0595165}}
\bptok{imsref}%
\end{bbook}
%
\endbibitem

\bibitem{Shaw1988}
%
\begin{barticle}[mr]
\bauthor{\bsnm{Shaw},~\bfnm{Daigee}\binits{D.}}
(\byear{1988}).
\btitle{On-site samples' regression: Problems of nonnegative integers,
truncation, and endogenous stratification}.
\bjournal{J. Econometrics}
\bvolume{37}
\bpages{211--223}.
\bid{doi={10.1016/0304-4076(88)90003-6}, issn={0304-4076}, mr={0932141}}
\bptok{imsref}%
\end{barticle}
%
\endbibitem

\bibitem{Skinner1994}
%
\begin{bincollection}[auto:STB|2011/10/17|13:52:43]
\bauthor{\bsnm{Skinner},~\bfnm{C.}\binits{C.}}
(\byear{1994}).
\btitle{Sample models and weights}.
In \bbooktitle{Proceedings of the Section on Survey Research Methods}
\bpages{133--142}.
\baddress{Washington, DC}: \bpublisher{American Statistical Association}.
\bptok{imsref}%
\end{bincollection}
%
\endbibitem

\bibitem{Sullivan2006}
%
\begin{bbook}[auto:STB|2011/10/17|13:52:43]
\bauthor{\bsnm{Sullivan},~\bfnm{P.}\binits{P.}},
\bauthor{\bsnm{Breidt},~\bfnm{F.}\binits{F.}},
\bauthor{\bsnm{Ditton},~\bfnm{R.}\binits{R.}},
\bauthor{\bsnm{Knuth},~\bfnm{B.}\binits{B.}},
\bauthor{\bsnm{Leaman},~\bfnm{B.}\binits{B.}},
\bauthor{\bsnm{O'Connell},~\bfnm{V.}\binits{V.}},
\bauthor{\bsnm{Parsons},~\bfnm{G.}\binits{G.}},
\bauthor{\bsnm{Pollock},~\bfnm{K.}\binits{K.}},
\bauthor{\bsnm{Smith},~\bfnm{S.}\binits{S.}} \AND
\bauthor{\bsnm{Stokes},~\bfnm{S.}\binits{S.}}
(\byear{2006}).
\btitle{Review of Recreational Fisheries Survey Methods}.
\baddress{Washington, DC}: \bpublisher{National Academies Press}.
\bptok{imsref}%
\end{bbook}
%
\endbibitem

\bibitem{vanderVaart1998}
%
\begin{bbook}[mr]
\bauthor{\bparticle{van~der} \bsnm{Vaart},~\bfnm{A.~W.}\binits{A.W.}}
(\byear{1998}).
\btitle{Asymptotic Statistics}.
\bseries{Cambridge Series in Statistical and Probabilistic Mathematics}
\bvolume{3}.
\baddress{Cambridge}: \bpublisher{Cambridge Univ. Press}.
\bid{mr={1652247}}
\bptok{imsref}%
\end{bbook}
%
\endbibitem

\end{thebibliography}
\end{document}